\documentclass[12pt,draft]{article} 
 \usepackage{ams}
    \font\twelvebm                       = cmmib10 at 12truept
    \font\tenbm                          = cmmib10 at 10truept
    \font\sevenbm                        = cmmib10 at 7truept
 = \tenbm \scriptscriptfont9              = \sevenbm
 
   at 10truept  at 10truept  at
 10truept  at 10truept
  at 10truept

 \mathchardef \BGamma            = "0900 \mathchardef \BDelta
 = "0901 \mathchardef \BTheta            = "0902 \mathchardef
 \BLambda           = "0903 \mathchardef \BXi               = "0904
 \mathchardef \BPi               = "0905 \mathchardef \BSigma
 = "0906 \mathchardef \BUpsilon          = "0907 \mathchardef \BPhi
 = "0908 \mathchardef \BPsi              = "0909 \mathchardef
 \BOmega            = "090A \mathchardef \Balpha            = "090B
 \mathchardef \Bbeta             = "090C \mathchardef \Bgamma
 = "090D \mathchardef \Bdelta            = "090E \mathchardef
 \Bepsilon          = "090F \mathchardef \Bzeta             = "0910
 \mathchardef \Beta              = "0911 \mathchardef \Btheta
 = "0912 \mathchardef \Biota             = "0913 \mathchardef
 \Bkappa            = "0914 \mathchardef \Blambda           = "0915
 \mathchardef \Bmu               = "0916 \mathchardef \Bnu
 = "0917 \mathchardef \Bxi               = "0918 \mathchardef \Bpi
 = "0919 \mathchardef \Brho              = "091A \mathchardef
 \Bsigma            = "091B \mathchardef \Btau              = "091C
 \mathchardef \Bupsilon          = "091D \mathchardef \Bphi
 = "091E \mathchardef \Bchi              = "091F \mathchardef \Bpsi
 = "0920 \mathchardef \Bomega            = "0921 \mathchardef
 \Bvarepsilon       = "0922 \mathchardef \Bvartheta         = "0923
 \mathchardef \Bvarpi            = "0924 \mathchardef \Bvarrho
 = "0925 \mathchardef \Bvarsigma         = "0926 \mathchardef
 \Bvarphi           = "0927
 \mathchardef \bA        = "0941 \mathchardef \bB        = "0942
 \mathchardef \bC        = "0943 \mathchardef \bD        = "0944
 \mathchardef \bE        = "0945 \mathchardef \bF        = "0946
 \mathchardef \bG        = "0947 \mathchardef \bH        = "0948
 \mathchardef \bI        = "0949 \mathchardef \bJ        = "094A
 \mathchardef \bK        = "094B \mathchardef \bL        = "094C
 \mathchardef \bM        = "094D \mathchardef \bN        = "094E
 \mathchardef \bO        = "094F \mathchardef \bP        = "0950
 \mathchardef \bQ        = "0951 \mathchardef \bR        = "0952
 \mathchardef \bS        = "0953 \mathchardef \bT        = "0954
 \mathchardef \bU        = "0955 \mathchardef \bV        = "0956
 \mathchardef \bW        = "0957 \mathchardef \bX        = "0958
 \mathchardef \bY        = "0959 \mathchardef \bZ        = "095A
 \mathchardef \ba        = "0961 \mathchardef \bb        = "0962
 \mathchardef \bc        = "0963 \mathchardef \bd        = "0964
 \mathchardef \bee       = "0965 
 \mathchardef \bff       = "0966 \mathchardef \bg        = "0967
 \mathchardef \bh        = "0968
 \mathchardef \bj        = "096A \mathchardef \bk        = "096B
 \mathchardef \bl        = "096C \mathchardef \bm        = "096D
 \mathchardef \bn        = "096E \mathchardef \bo        = "096F
 \mathchardef \bp        = "0970 \mathchardef \bq        = "0971
 \mathchardef \br        = "0972 \mathchardef \bs        = "0973
 \mathchardef \bt        = "0974 \mathchardef \bu        = "0975
 \mathchardef \bv        = "0976 \mathchardef \bw        = "0977
 \mathchardef \bx        = "0978 \mathchardef \by        = "0979
 \mathchardef \bz        = "097A

 \font\tencb            = cmssbx10 scaled \magstep4 \font\eigcb
 = cmssbx10 scaled \magstep2 \textfont8             = \tencb
 \mathchardef\bAs       = "1841
 \def\Asem#1#2{\mathop{\vrule height10.5pt depth5.5pt width0pt\bAs}_{#1}^{#2}}
 \def\asem#1#2{
          \ifmmode
          \ifinner
            \raise0.9pt\hbox{$\scriptstyle\bAs$}_{#1}^{#2}
          \else
            \Asem{#1}{#2}
         \fi
          \fi
          }

 \newtheorem{theo}{\small\bf Theorem}[section]
 \newtheorem{lem}{\small\bf Lemma}[section]
 \newtheorem{prop}{\small\bf Proposition}[section]
 \newtheorem{rem}{\small\bf Remark}[section]
 \newenvironment{REM}{\begin{rem} \rm}{\end{rem}}
 \newtheorem{defi}{\small\bf Definition}[section]
 \newenvironment{DEFI}{\begin{defi} \rm}{\end{defi}}
 \newtheorem{cor}{\small\bf Corollary}[section]
 
 \renewcommand{\Pr}{\mbox{\rm\hspace*{.2ex}I\hspace{-.5ex}P\hspace*{.2ex}}}

 \newcommand{\be}{\begin{equation}}
 \newcommand{\ee}{\end{equation}}
 \newcommand{\law}{\stackrel{\mbox{\footnotesize d}}{=}}

 \newcommand{\E}{\mbox{\rm \hspace*{.2ex}I\hspace{-.5ex}E\hspace*{.2ex}}}
 \newcommand{\Var}{\mbox{\rm \hspace*{.2ex}Var\hspace*{.2ex}}}
 \newcommand{\Cov}{\mbox{\rm \hspace*{.2ex}Cov\hspace*{.2ex}}}
 \newcommand{\MSE}{\mbox{\rm \hspace*{.2ex}MSE\hspace*{.2ex}}}
 \newcommand{\D}{\mbox{\rm \hspace*{.2ex}I\hspace{-.5ex}D\hspace*{.2ex}}}
 
 \newcommand{\X}{X_1,X_2,\ldots ,X_n}
 \newcommand{\Xo}{X_{1:n},X_{2:n},\ldots ,X_{n:n}}
 
 \newcommand{\Xp}{X_{1:1},X_{2:2},\ldots ,X_{n:n}}
 \newcommand{\Xs}{X^*_1,X^*_2,\ldots ,X^*_n}
 \newcommand{\Xso}{X^*_{1:n},X^*_{2:n},\ldots ,X^*_{n:n}}
 \newcommand{\Xsoo}{X^*_{1:n}\leq X^*_{2:n}\leq \cdots \leq X^*_{n:n}}
 \newcommand{\Xsp}{X^*_{1:1},X^*_{2:2},\ldots ,X^*_{n:n}}
 \newcommand{\Xspo}{X^*_{1:1}\leq X^*_{2:2}\leq \cdots \leq X^*_{n:n}}
 \newcommand{\Xpo}{X_{1:1}\leq X_{2:2}\leq \cdots \leq X_{n:n}}
 \newcommand{\Xspmo}{X^*_{1:1}\geq X^*_{1:2}\geq \cdots \geq X^*_{1:n}}

 \newcommand{\bbb}[1]{\mbox{\boldmath $ #1 $}}

 \newenvironment{pr}[1]{{\small\bf {#1}:}}{}

 \title{ \Large\bf Linear Estimation of Location
 and Scale Parameters \\
 using Partial Maxima
 \vspace*{-.6em}}

 \author{\large
 Nickos Papadatos}

 \date{\small
 Department of Mathematics, Section of Statistics and O.R.,
 University of Athens, \\ Panepistemiopolis, 157 84 Athens, Greece.
 {\tt e-mail:\ npapadat@math.uoa.gr}\vspace*{-1em}}

\begin{document}

 \maketitle
 \vspace*{-2em}

 \thispagestyle{empty}

 \begin{abstract}
  \noindent
  Consider an i.i.d.\ sample $\Xs$ from a location-scale family, and
  assume that the only available observations consist of the
  partial maxima (or minima)
  sequence, $\Xsp$, where
  $X^*_{j:j}=\max\{ X^*_1,  \ldots,X^*_j \}$.
  This kind of truncation
  appears in several circumstances, including
  best performances in athletics events.
  In the case of partial maxima, the form of the BLUEs (best linear
  unbiased estimators)
  is quite similar to the form of the well-known Lloyd's (1952, Least-squares estimation
  of location and scale parameters using order statistics, {\it Biometrika}, vol.\ 39, pp.\ 88--95) BLUEs,
  based on (the sufficient sample of) order statistics,
  but, in contrast to the
  classical case,
  their consistency is no longer obvious.
  The present paper is mainly concerned with
  the scale
  parameter,
  showing
  that the
  variance of the partial maxima BLUE
  is at most
  of order $O(1/\log n)$, for a wide class of distributions.
  \end{abstract}
 \vspace*{-1em}

 \section{Introduction}\vspace*{-.5em}
 \label{sec1}

 There are several situations
 where the ordered random
 sample,
 \be \label{1.1}
 \Xsoo,
 \ee
 corresponding to the i.i.d.\
 random sample, $\Xs$, is not fully reported, because the
 values of interest are the higher (or lower), up-to-the-present,
 record values based on the initial sample, i.e., the partial
 maxima (or minima) sequence
 \be \label{1.2}
 \Xspo,
 \ee
 where $X^*_{j:j}=\max\{ X^*_1,\ldots,X^*_j \}$.
 A situation of this kind commonly appears in athletics, when only
 the best performances are recorded.

 Through this article we assume that the i.i.d.\ data arise from a
 location-scale family,
 \[
 \{ F((\cdot-\theta_1)/\theta_2); \, \theta_1\in\R, \theta_2 >0 \},
 \]
 where the d.f.\ $F(\cdot)$ is free of parameters and has finite,
 non-zero variance (so that $F$ is non-degenerate), and we
 consider the partial maxima BLUE (best linear unbiased
 estimator) for both parameters $\theta_1$ and $\theta_2$. This
 consideration is along the lines of the classical Lloyd's (1952)
 BLUEs, the only difference being that the linear estimators are
 now based on the ``insufficient sample'' (\ref{1.2}), rather than
 (\ref{1.1}), and this fact implies a substantial reduction on the
 available information. Tryfos and Blackmore (1985) used
 this kind of data to predict future records in athletic events,
 Samaniego and Whitaker (1986, 1988) estimated the population characteristics,
 while Hofmann and Nagaraja (2003) investigated the amount of
 Fisher Information contained in such data; see also
 Arnold, Balakrisnan \& Nagaraja (1998, Section 5.9).

 A natural question concerns the consistency of the resulting
 BLUEs, since too much lack of information presumably would result
 to inconsistency (see at the end of Section 6).
 Thus, our main focus is on
 conditions guaranteeing consistency,
 and the main result
 shows that this is indeed the case for the
 scale
 parameter BLUE from a wide class of distributions.
 Specifically, it is
 shown that the variance of the BLUE
 is at most of order $O(1/\log n)$, when
 $F(x)$ has a log-concave
 density $f(x)$
 and satisfies the Von Mises-type condition (\ref{5.11}) or
 (\ref{6.1}) (cf.\ Galambos (1978)) on the right end-point of its
 support (Theorem \ref{theo5.2}, Corollary \ref{cor6.1}). The
 result is applicable to several commonly used distributions, like
 the Power distribution (Uniform), the Weibull (Exponential), the
 Pareto, the Negative Exponential, the Logistic,
 the Extreme Value (Gumbel)
 and the Normal (see
 section \ref{sec6}).
 A consistency result for the partial maxima BLUE of
 the location parameter
 would be desirable to be
 included here, but it seems that the proposed technique (based on
 partial maxima spacings, section \ref{sec4}) does not suffice for
 deriving it. Therefore, the consistency for the location parameter
 remains an open problem in general, and it is just highlighted by
 a particular application to the Uniform location-scale family
 (section \ref{sec3}).

 The proof of the main result
 depends on the fact that, under mild conditions, the partial
 maxima spacings have non-positive correlation. The class of distributions
 having this
 property is called NCP (negative correlation for partial maxima
 spacings).
 It is shown here that any log-concave distribution
 with finite variance belongs to NCP
 (Theorem
 \ref{theo4.2}).
 In particular, if a distribution function
 has a density which is
 either log-concave or non-increasing
 then it is a member of NCP.
 For ordinary spacings, similar sufficient conditions
 were shown by
 Sarkadi (1985) and
 Bai, Sarkar \& Wang (1997) --
 see also David and Nagaraja (2003, pp.\ 187--188), Burkschat (2009),
 Theorem 3.5 --
 and will be referred as
 ``S/BSW-type conditions''.

 In every experiment where the i.i.d.\ observations arise in a sequential manner,
 the partial maxima data describe the
 best performances in a natural way,
 as the experiment goes on,
 in contrast to the first $n$
 record values, $R_1,R_2,\ldots,R_n$, which are obtained
 from an inverse sampling scheme -- see, e.g., Berger and Gulati (2001).
 Due to the very rare
 appearance of records, in the latter case it is implicitly assumed
 that the sample size is, roughly, $e^n$.
 This has a similar effect in the partial maxima setup, since
 the number of different values
 are about $\log n$, for large sample size $n$.
 Clearly, the total amount of information in the partial
 maxima sample is the same as that given by the (few)
 record values augmented by record times.
 The essential difference of these models (records / partial
 maxima) in statistical applications
 is highlighted, e.g., in Tryfos and Blackmore (1985),
 Samaniego and Whitaker (1986, 1988), Smith (1988), Berger and Gulati (2001) and
 Hofmann and Nagaraja
 (2003) -- see also Arnold,Balakrishnan \& Nagaraja
 (1998, Chapter 5).

 \section{Linear estimators based on partial maxima}\vspace*{-.5em}
 \setcounter{equation}{0} \label{sec2}

 Consider the random
 sample $\Xs$ from $F((x-\theta_1)/\theta_2)$  and the corresponding
 partial maxima sample $\Xspo$
 ($\theta_1\in\R$ is
 the location parameter and
 $\theta_2>0$ is the scale parameter; both parameters are unknown).
 Let also $\X$ and $\Xpo$ be the
 corresponding
 samples from the completely specified d.f.\ $F(x)$,
 that generates the location-scale family. Since
 \[
 ( \Xsp )' \law ( \theta_1+\theta_2 X_{1:1}, \theta_1+\theta_2 X_{2:2},
  \ldots ,\theta_1+\theta_2 X_{n:n} )',
 \]
 a linear estimator based  on
 partial maxima has the form
 \[
 L=\sum_{i=1}^{n} c_i X^*_{i:i} \law \theta_1 \sum_{i=1}^{n} c_i
 + \theta_2 \sum_{i=1}^{n} c_i X_{i:i},
 \]
 for some constants $c_i$, $i=1,2,\ldots,n$.

 Let $\bbb{X}=(\Xp)'$ be the random vector of partial maxima from
 the known d.f.\ $F(x)$, and use the notation
 \be \label{2.1}
 \bbb{\mu}=\E [\bbb{X}], \hspace*{2ex} \bbb{\Sigma}= \D[\bbb{X}]
 \hspace*{1ex}\mbox{and}\hspace*{1ex} \bbb{E}=\E[ \bbb{X}
 \bbb{X}'],
 \ee
 where $\D[\bbb{\xi}]$ denotes the dispersion matrix
 of any random vector $\bbb{\xi}$. Clearly,
 \[
 \bbb{\Sigma}=\bbb{E}-\bbb{\mu} \bbb{\mu}',\hspace*{2ex} \bbb{\Sigma}>0,
 \hspace*{2ex} \bbb{E}>0.
 \]
 The linear estimator $L$
 is called BLUE for $\theta_k$
 ($k=1,2$) if it is unbiased for $\theta_k$ and its variance is
 minimal, while it is called BLIE (best linear invariant estimator)
 for $\theta_k$ if it is invariant for $\theta_k$ and its mean
 squared error, $\MSE[L]=\E[L-\theta_k]^2$, is minimal.
 Here ``invariance" is understood in the sense of location-scale invariance as it is defined, e.g., in Shao (2005, p.\ xix).

 Using the above notation it is easy to verify the following formulae for the BLUEs
 and their variances. They are the partial maxima analogues of Lloyd's (1952) estimators
 and, in the case of partial minima, have been obtained
 by Tryfos and Blackmore (1985), using least squares. A proof is attached
 here for easy reference.
 \begin{prop}
 \label{prop2.1}
 The partial maxima {\rm BLUEs} for $\theta_1$ and for $\theta_2$ are,
 respectively,
  \be \label{2.2}
   L_1=-\frac{1}{\Delta}\bbb{\mu}'\bbb{\Gamma}\bbb{X}^*\ \ \mbox{and} \ \
   L_2=\frac{1}{\Delta}{\bf 1}'\bbb{\Gamma}\bbb{X}^*,
  \ee
 where $\bbb{X}^*=(\Xsp )'$,
 $\Delta=({\bf 1}'\bbb{\Sigma}^{-1}{\bf 1}) (\bbb{\mu}'\bbb{\Sigma}^{-1}
 \bbb{\mu})-({\bf 1}'\bbb{\Sigma}^{-1}\bbb{\mu})^2>0$,
 ${\bf 1}=(1,1,\ldots,1)'\in \R^n$ and
 $\bbb{\Gamma}=\bbb{\Sigma}^{-1}({\bf 1}\bbb{\mu}' -\bbb{\mu}{\bf
 1}')\bbb{\Sigma}^{-1}$. The corresponding  variances are
 \be
 \label{2.3}
 \Var [L_1] =
 \frac{1}{\Delta}(\bbb{\mu}'\bbb{\Sigma}^{-1}\bbb{\mu}) \theta_2^2
 \ \ \mbox{and}\ \
  \Var [L_2] =
 \frac{1}{\Delta}({\bf 1}'\bbb{\Sigma}^{-1}{\bf 1}) \theta_2^2.
 \ee
 \end{prop}
 \begin{pr}{Proof}
 Let $\bbb{c}=(c_1,c_2,\ldots,c_n)'\in \R^n$ and
 $L=\bbb{c}'\bbb{X}^*$. Since
 $\E[L]=(\bbb{c}'\bbb{1})\theta_1+(\bbb{c}'\bbb{\mu})\theta_2$, $L$
 is unbiased for $\theta_1$ iff $\bbb{c}'{\bf 1}=1$ and
 $\bbb{c}'\bbb{\mu}=0$, while it is unbiased for $\theta_2$ iff
 $\bbb{c}'{\bf 1}=0$ and $\bbb{c}'\bbb{\mu}=1$. Since
 $\Var[L]=(\bbb{c}'\bbb{\Sigma}\bbb{c})\theta_2^2$, a simple
 minimization argument for $\bbb{c}'\bbb{\Sigma}\bbb{c}$ with
 respect to $\bbb{c}$, using Lagrange multipliers, yields the
 expressions (\ref{2.2}) and (\ref{2.3}). $\Box$
 \end{pr}
 \bigskip

 Similarly, one can derive the partial maxima version of Mann's (1969)
 best linear invariant estimators (BLIEs),  as follows.

 \begin{prop}
 \label{prop2.2}
 The partial maxima {\rm BLIEs} for $\theta_1$ and for $\theta_2$ are,
 respectively,
  \be \label{2.4}
   T_1=\frac{{\bf 1}'\bbb{E}^{-1}\bbb{X}^*}{{\bf 1}'\bbb{E}^{-1}{\bf 1}}
   \ \ \mbox{and} \ \
   T_2=\frac{{\bf 1}'\bbb{G}\bbb{X}^*}{{\bf 1}'\bbb{E}^{-1}{\bf 1}},
  \ee
 where $\bbb{X}^*$ and ${\bf 1}$ are as in Proposition {\rm
 \ref{prop2.1}} and
 $\bbb{G}=\bbb{E}^{-1}({\bf 1}\bbb{\mu}'-\bbb{\mu}{\bf 1}')\bbb{E}^{-1}$.
 The corresponding  mean squared errors are
 \be \label{2.5}
 \MSE
 [T_1] = \frac{\theta_2^2}{{\bf 1}'\bbb{E}^{-1}{\bf 1}} \ \
 \mbox{and}\ \
 \MSE [T_2] =
 \left( 1-\frac{D}{{\bf 1}'\bbb{E}^{-1}{\bf 1}}\right) \theta_2^2,
 \ee
 where $D=({\bf 1}'\bbb{E}^{-1}{\bf 1}) (\bbb{\mu}'\bbb{E}^{-1}
 \bbb{\mu})-({\bf 1}'\bbb{E}^{-1}\bbb{\mu})^2>0$.
 \end{prop}
 \begin{pr}{Proof}
 Let $L=L(\bbb{X}^*)=\bbb{c}'\bbb{X}^*$ be an arbitrary linear
 statistic. Since 
 $L(b\bbb{X}^*+a{\bf 1})=a (\bbb{c}'{\bf 1})+b
 L(\bbb{X}^*)$ for arbitrary $a\in\R$ and $b>0$, it follows that $L$ is invariant for $\theta_1$ 
 iff $\bbb{c}'{\bf 1}=1$ while it is invariant for $\theta_2$ iff $\bbb{c}'{\bf 1}=0$. 
 Both (\ref{2.4}) and
 (\ref{2.5}) now follow by a simple minimization argument, since in
 the first case we have to minimize the mean squared error
 $\E[L-\theta_1]^2=(\bbb{c}'\bbb{E}\bbb{c})\theta_2^2$ under
 $\bbb{c}'{\bf 1}=1$, while in the second one, we have to minimize
 the mean squared error $\E[L-\theta_2]^2=
 (\bbb{c}'\bbb{E}\bbb{c}-2\bbb{\mu}'\bbb{c}+1)\theta_2^2$ under
 $\bbb{c}'{\bf 1}=0$.
 $\Box$
 \end{pr}
 \bigskip

 The above formulae (\ref{2.2})-(\ref{2.5}) are well-known for
 order statistics and records -- see
 David (1981, Chapter 6),
 Arnold, Balakrishnan \&
 Nagaraja (1992, Chapter 7; 1998, Chapter 5),
 David and Nagaraja (2003, Chapter 8). In the
 present setup, however, the meaning of $\bbb{X}^*$, $\bbb{X}$,
 $\bbb{\mu}$, $\bbb{\Sigma}$ and $\bbb{E}$ is completely different.
 In the case of order statistics, for example, the vector
 $\bbb{\mu}$, which is the mean vector of the order statistics
 $\bbb{X}=(\Xo)'$ from the known distribution $F(x)$, depends on the
 sample size $n$, in the sense that the components of the vector
 $\bbb{\mu}$ completely change with $n$. In the present case of
 partial maxima, the first $n$ entries of the vector $\bbb{\mu}$,
 which is the mean vector of the partial maxima $\bbb{X}=(\Xp)'$ from
 the known distribution $F(x)$, remain constant for all sample sizes
 $n'$ greater than or equal to $n$. Similar observations apply for
 the matrices $\bbb{\Sigma}$ and $\bbb{E}$. This fact seems to be
 quite helpful for the construction of tables giving the means,
 variances and covariances of partial maxima for samples up to a size
 $n$. It should be noted, however, that even when $F(x)$ is
 absolutely continuous with density $f(x)$ (as is usually the case
 for location-scale families), the joint distribution of $(X_{i:i},
 X_{j:j})$ has a singular part, since
 $\Pr[X_{i:i}=X_{j:j}]=i/j>0$, $i<j$.
 Nevertheless, there exist simple
 expectation and covariance
 formulae (Lemma \ref{lem2.2}).

 As in the order statistics setup, the actual application of
 formulae (\ref{2.2}) and (\ref{2.4}) requires closed forms for
 $\bbb{\mu}$ and $\bbb{\Sigma}$, and also to invert the $n\times n$
 matrix $\bbb{\Sigma}$. This can be done only for very particular
 distributions (see next section, where we apply the results to the
 Uniform distribution). Therefore, numerical methods should be
 applied in general. This, however, has a theoretical cost: It
 is not a trivial fact to verify consistency of the estimators, even in the
 classical case of order statistics.
 The main purpose of this article
 is in verifying consistency for the partial maxima BLUEs.
 Surprisingly, it seems that a solution of this problem is not
 well-known, at least to our knowledge, even for the classical BLUEs
 based on order statistics. However, even if the result of the
 following lemma is known, its proof has an independent interest,
 because it proposes alternative (to BLUEs)
 $n^{-1/2}$--consistent unbiased linear estimators
 and provides the intuition for the derivation of the main
 result of the present article.

 \begin{lem}\label{lem2.1}
 The classical {\rm BLUEs}  of $\theta_1$ and $\theta_2$,
 based on order statistics from a location-scale family,
 created by a distribution $F(x)$
 with finite non-zero variance, are consistent. Moreover,
 their variance is at most of order $O(1/n)$.
 \end{lem}
 \begin{pr}{Proof}
 Let $\bbb{X}^*=(\Xso)'$ and $\bbb{X}=(\Xo)'$ be the ordered samples
 from $F((x-\theta_1)/\theta_2)$ and $F(x)$, respectively, so
 that $\bbb{X}^*\law \theta_1 {\bf 1}+\theta_2 \bbb{X}$.
 Also write $X_1^*,X_2^*,\ldots,X_n^*$ and
 $X_1,X_2,\ldots,X_n$ for the corresponding i.i.d.\ samples.
 We consider the linear estimators
 \[
 S_1=\overline{X}^*=\frac1n\sum_{i=1}^n X_i^*\law \theta_1+
 \theta_2 \overline{X}
 \]
 and
 \[
 S_2=\frac{1}{n(n-1)}\sum_{i=1}^n\sum_{j=1}^n |X^*_j-X^*_i|\law
 \frac{\theta_2}{n(n-1)}\sum_{i=1}^n\sum_{j=1}^n |X_j-X_i|,
 \]
 i.e., $S_1$ is the sample mean and $S_2$ is a multiple of
 Gini's statistic. Observe that
 both $S_1$ and $S_2$ are linear estimators in order statistics.
 [In particular, $S_2$ can be written as
 $S_2=4(n(n-1))^{-1}\sum_{i=1}^n (i-(n+1)/2)X_{i:n}^*$.] 
 Clearly,
 $\E(S_1)=\theta_1+\theta_2\mu_0$, $\E(S_2)=\theta_2 \tau_0$,
 where $\mu_0$ is the mean, $\E(X_1)$, of the distribution $F(x)$
 and $\tau_0$ is the positive finite parameter
 $\E|X_1-X_2|$. Since $F$ is known, both $\mu_0\in\R$ and $\tau_0>0$
 are known constants,
 and we can construct the linear estimators
 $U_1=S_1-(\mu_0/\tau_0)S_2$ and $U_2=S_2/\tau_0$. Obviously,
 $\E(U_k)=\theta_k$, $k=1,2$, and both $U_1$, $U_2$ are linear
 estimators of the form $T_n=(1/n)\sum_{i=1}^n \delta(i,n) X^*_{i:n}$,
 with $|\delta(i,n)|$ uniformly bounded for all $i$ and $n$.
 If $\sigma_0^2$ is the (assumed finite) variance of $F(x)$,
 it follows that
 \begin{eqnarray*}
 \Var[T_n] & \leq & \displaystyle \frac{1}{n^2} \sum_{i=1}^n \sum_{j=1}^n
 |\delta(i,n)||\delta(j,n)| \Cov(X^*_{i:n},X^*_{j:n}) \\
 & \leq & \displaystyle \frac{1}{n^2}
 \left(\max_{1\leq i \leq n} |\delta(i,n)|\right)^2
 \Var(X^*_{1:n}+X^*_{2:n}+\cdots+X^*_{n:n}) \\
 & = & \displaystyle \frac{1}{n}
 \left(\max_{1\leq i \leq n} |\delta(i,n)|\right)^2
 \theta_2^2
 \sigma_0^2=O(n^{-1})\to 0,\ \ \ \mbox{as}\ \ n\to\infty,
 \end{eqnarray*}
 showing that $\Var(U_k)\to 0$, and thus $U_k$ is consistent
 for $\theta_k$, $k=1,2$. Since
 $L_k$ has minimum variance among all linear unbiased estimators,
 it follows that $\Var(L_k)\leq \Var (U_k)\leq O(1/n)$,
 and the result follows. $\Box$
 \end{pr}
 \bigskip

 The above lemma implies that the mean
 squared error of the BLIEs, based on order statistics, is at most
 of order $O(1/n)$, since they have smaller mean squared error than
 the BLUEs, and thus they are also consistent. More important is
 the fact that, with the technique used in Lemma \ref{lem2.1}, one
 can avoid all computations involving means, variances and
 covariances of order statistics, and it does not need to invert
 any matrix, in order to prove consistency (and in order to obtain
 $O(n^{-1})$-consistent estimators). Arguments of similar
 kind will be applied in section \ref{sec5}, when the problem of
 consistency for the partial maxima BLUE of $\theta_2$ will be
 taken under consideration.

 We now turn in the partial maxima case.
 Since actual application of partial maxima BLUEs and BLIEs
 requires the computation of the first two moments
 of $\bbb{X}=(\Xp)'$ in terms of the
 completely specified d.f.\ $F(x)$, the following formulae
 are to be mentioned here (cf.\ Jones and Balakrishnan (2002)).

 \begin{lem} \label{lem2.2} Let $\Xpo$ be the partial maxima sequence
 based on an arbitrary d.f.\ $F(x)$. \\
 {\rm (i)} For $i\leq j$, the joint d.f.\ of $(X_{i:i},X_{j:j})$ is
 \be\label{2.6}
 F_{X_{i:i},X_{j:j}}(x,y)= \left\{
 \begin{array}{lll}
 F^j(y) & \mbox { if } & x\geq y, \\
 F^i(x)F^{j-i}(y) & \mbox { if } & x\leq y.
 \end{array}
 \right. \ee {\rm (ii)} If $F$ has finite first moment, then
 \be
 \label{2.7}
 \mu_i=\E[X_{i:i}]=\int_0^{\infty} (1-F^i(x))\ dx -
 \int_{-\infty}^0 F^i(x) dx
 \ee
 is finite for all $i$. \\
 {\rm (iii)} If $F$ has finite second moment, then
 \be\label{2.8}
 \sigma_{ij}=\Cov[X_{i:i},X_{j:j}]= \int\int_{-\infty<x<y<\infty}
 F^i(x) (F^{j-i}(x)+F^{j-i}(y)) (1-F^i (y)) \ dy \ dx
 \ee is finite
 and non-negative for all $i\leq j$. In particular,
 \be\label{2.9}
 \sigma_{ii}=\sigma_i^2=\Var[X_{i:i}]= 2
 \int\int_{-\infty<x<y<\infty} F^i(x) (1-F^i (y)) \ dy \ dx.
 \ee
 \end{lem}
 \begin{pr}{Proof} (i) is trivial and (ii) is well-known. (iii)
 follows from Hoeffding's identity
 \[
 \Cov[X,Y]=\int_{-\infty}^{\infty}\int_{-\infty}^{\infty}
 (F_{X,Y}(x,y)-F_X(x)F_Y(y)) \ dy \ dx
 \]
 (see Hoeffding (1940), Lehmann (1966), Jones and Balakrishnan
 (2002), among others), applied to $(X,Y)=(X_{i:i},X_{j:j})$ with
 joint d.f.\ given by (\ref{2.6}) and marginals $F^i(x)$ and
 $F^j(y)$. $\Box$
 \end{pr}
 \bigskip

 Formulae (\ref{2.7})-(\ref{2.9}) enable the computation
 of means, variances and covariances of partial maxima, even in the
 case where the distribution $F$ does not have a density. Tryfos and
 Blackmore (1985) obtained an expression for the covariance of partial
 minima involving means and covariances of order statistics from
 lower sample sizes.

 \section{A tractable case: the Uniform location-scale family}\vspace*{-.5em}
 \setcounter{equation}{0}\label{sec3}

 Let $\Xs \sim U(\theta_1,\theta_1+\theta_2)$, so that $(\Xsp)'\law
 \theta_1 {\bf 1} + \theta_2 \bbb{X}$, where $\bbb{X}=(\Xp)'$ is
 the partial maxima sample from the standard Uniform distribution.
 Simple calculations, using (\ref{2.7})-(\ref{2.9}), show that the
 mean vector $\bbb{\mu}=(\mu_i)$ and the dispersion matrix
 $\bbb{\Sigma}=(\sigma_{ij})$ of $\bbb{X}$ are given by
 (see also Tryfos and Blackmore (1985), eq.\ (3.1))
 \[
 \mu_i=\frac{i}{i+1}\ \ \mbox{and}\ \ \sigma_{ij}=\frac{i}{(i+1)(j+1)(j+2)}
 \ \ \mbox{for}\ \  1\leq i\leq j\leq n.
 \]
 Therefore, $\bbb{\Sigma}$ is a patterned matrix of the
 form $\sigma_{ij}=a_i b_j$
 for $i\leq j$, and thus, its inverse is tridiagonal; see Graybill
 (1969, Chapter 8), Arnold, Balakrishnan \& Nagaraja (1992,
 Lemma 7.5.1). Specifically,
 \[
 \bbb{\Sigma}^{-1}=\left(
 \begin{array}{cccccc}
 \gamma_1 & -\delta_1 & 0 & \ldots & 0 & 0\\
 -\delta_1 & \gamma_2 & -\delta_2 & \ldots & 0 & 0\\
 0 & -\delta_2 & \gamma_3 & \ldots & 0 & 0\\
 \vdots & \vdots& \vdots & \vdots & \vdots & \vdots \\
 0 & 0 & 0 & \ldots & \gamma_{n-1} & -\delta_{n-1} \\
 0 & 0 & 0 & \ldots & -\delta_{n-1} & \gamma_n
 \end{array}
 \right)
 \]
 where
 \begin{eqnarray*}
 \displaystyle \gamma_i=\frac{4(i+1)^3(i+2)^2}{(2i+1)(2i+3)}, \ \
 \delta_i =\frac{(i+1)(i+2)^2(i+3)}{2i+3}, \ \ i=1,2,\ldots,n-1, &
 \\
 \mbox{ and } \ \
  \displaystyle  \gamma_n=\frac{(n+1)^2 (n+2)^2}{2n+1}. \hspace*{30ex} &
 \end{eqnarray*}
 Setting $a(n)={\bf 1}' \bbb{\Sigma}^{-1}{\bf 1}$,
 $b(n)=({\bf 1}-\bbb{\mu})' \bbb{\Sigma}^{-1}({\bf 1}-\bbb{\mu})$
 and $c(n)=({\bf 1}-\bbb{\mu})' \bbb{\Sigma}^{-1}{\bf 1}$, we get

 \begin{eqnarray*}
 a(n) & = & \frac{(n+1)^2 (n+2)^2}{2n+1}-2\sum_{i=1}^{n-1}
 \frac{(i+1)(i+2)^2(3i+1)}{(2i+1)(2i+3)}=n^2+o(n^2), \\
 b(n) & = & \frac{(n+2)^2}{2n+1}-2\sum_{i=1}^{n-1}
 \frac{(i-1)(i+2)}{(2i+1)(2i+3)}=\frac12 \log n +o(\log n), \\
 c(n) & = & \frac{(n+1) (n+2)^2}{2n+1}-\sum_{i=1}^{n-1}
 \frac{(i+2)(4i^2+7i+1)}{(2i+1)(2i+3)}=n+o(n).
 \end{eqnarray*}
 Applying (\ref{2.3}) we obtain
 \begin{eqnarray*}
 \Var[L_1] & = & \frac{a(n)+b(n)-2c(n)}{a(n)b(n)-c^2(n)}\theta_2^2
 =\left(\frac{2}{\log n}+o\left( \frac{1}{\log n}\right)\right)\theta_2^2,
 \ \ \mbox{and}\ \ \\
 \Var[L_2]& = & \frac{a(n)}{a(n)b(n)-c^2(n)}\theta_2^2
 =\left(\frac{2}{\log n}+o\left( \frac{1}{\log n}\right)\right)\theta_2^2.
 \end{eqnarray*}

 The preceding computation shows that, for the Uniform location-scale family,
 the partial maxima BLUEs are consistent for both the location and
 the scale parameters, since their variance goes to zero with the
 speed of $2/\log n$. This fact, as expected, contradicts the
 behavior of the ordinary order statistics BLUEs, where the
 speed of convergence is of order $n^{-2}$ for the variance
 of both Lloyd's estimators. However, the comparison is
 quite unfair here, since Lloyd's estimators are based on the
 complete sufficient statistic $(X^*_{1:n},X^*_{n:n})$, and thus
 the variance of order statistics BLUE is minimal among all
 unbiased estimators.

 On the other hand we should emphasize that, under the same model,
 the BLUEs (and the BLIEs) based solely on the first $n$ upper records
 are not even consistent. In fact, the variance of both BLUEs converges to
 $\theta^2_2/3$, and the MSE of both BLIEs approaches $\theta_2^2/4$,
 as $n\to\infty$; see Arnold, Balakrishnan \& Nagaraja (1998,
 Examples 5.3.7 and 5.4.3).

 \section{Scale estimation and partial maxima spacings}\vspace*{-.5em}
 \setcounter{equation}{0}\label{sec4}

 In the classical order statistics setup,
 Balakrishnan and Papadatos (2002) observed
 that the computation of BLUE (and
 BLIE) of the scale parameter is simplified considerably
 if one uses spacings instead of order statistics -- cf.\
 Sarkadi (1985).
 Their observation applies here too,
 and simplifies the form of the partial maxima BLUE (and BLIE).

 Specifically, define the partial maxima spacings
 as $Z_i^*=X^*_{i+1:i+1}-X^*_{i:i}\geq 0$ and
 $Z_i=X_{i+1:i+1}-X_{i:i}\geq 0$,
 for $i=1,2,\ldots,n-1$, and let
 $\bbb{Z}^*=(Z_1^*,Z_2^*,\ldots,Z_{n-1}^*)'$ and
 $\bbb{Z}=(Z_1,Z_2,\ldots,Z_{n-1})'$. Clearly,
 $\bbb{Z}^* \law \theta_2 \bbb{Z}$, and any unbiased (or even
 invariant) linear estimator of $\theta_2$ based on
 the partial maxima sample, $L=\bbb{c}'\bbb{X}^*$,
 should necessarily satisfy $\sum_{i=1}^n c_i=0$
 (see the proofs of Propositions \ref{prop2.1} and \ref{prop2.2}).
 Therefore, $L$ can be expressed as a linear function on $Z^*_i$'s,
 $L=\bbb{b}'\bbb{Z}^*$, where now
 $\bbb{b}=(b_1,b_2,\ldots,b_{n-1})'\in \R^{n-1}$. Consider
 the mean vector
 $\bbb{m}=\E[\bbb{Z}]$, the dispersion matrix
 $\bbb{S}=\D[\bbb{Z}]$, and the second moment matrix
 $\bbb{D}=\E[\bbb{Z}\bbb{Z}']$ of $\bbb{Z}$. Clearly,
 $\bbb{S}=\bbb{D}-\bbb{m}\bbb{m}'$, $\bbb{S}>0$, $\bbb{D}>0$, and
 the vector $\bbb{m}$ and the matrices $\bbb{S}$ and $\bbb{D}$ are
 of order $n-1$. Using exactly the same arguments as in
 Balakrishnan and Papadatos (2002), it is easy to verify the following.

 \begin{prop}\label{prop4.1}
 The partial maxima {\rm BLUE} of $\theta_2$, given in Proposition
 {\rm \ref{prop2.1}}, has the alternative form \vspace*{-0.6em}
 \be\label{4.1}
 L_2=\frac{\bbb{m}'\bbb{S}^{-1}\bbb{Z}^*}{\bbb{m}'\bbb{S}^{-1}\bbb{m}},
 \  \  \mbox{ with }\ \
 \Var[L_2]=\frac{\theta_2^2}{\bbb{m}'\bbb{S}^{-1}\bbb{m}}, \ee
 while the corresponding {\rm BLIE}, given in Proposition {\rm
 \ref{prop2.2}}, has the alternative form
 \be\label{4.2}
 T_2=\bbb{m}'\bbb{D}^{-1}{\bbb{Z}}^*, \  \  \mbox{
 with }\ \ \MSE[T_2]=(1-\bbb{m}'\bbb{D}^{-1}\bbb{m})\theta_2^2.
 \ee
 \end{prop}


 It should be noted that, in general, the non-negativity of the
 BLUE of $\theta_2$ does not follow automatically,
 even for order
 statistics. In the order statistics setup, this problem
 was posed  by Arnold, Balakrishnan \&
 Nagaraja (1992), and the best known
 result, till now, is the one
 given by Bai, Sarkar \& Wang (1997) and Sarkadi (1985).
 Even after the slight
 improvement, given by Balakrishnan and Papadatos (2002) and
 by Burkschat (2009),
 the general case remains unsolved. The same
 question (of non-negativity of the BLUE) arises in the partial
 maxima setup, and the following theorem provides a partial
 positive answer.
 We omit the
 proof, since it again
 follows by a straightforward
 application of the arguments given in Balakrishnan and
 Papadatos
 (2002). \vspace*{-0.6em}

 \begin{theo}
 \label{theo4.1} {\rm (i)} There exists a constant $a=a_n(F)$,
 $0<a<1$, depending only on the sample size $n$ and the d.f.\
 $F(x)$ {\rm (i.e., $a$ is free of the parameters $\theta_1$ and
 $\theta_2$)}, such that $ T_2=a \, L_2. $ This constant is given
 by $ a=\bbb{m}'\bbb{D}^{-1}\bbb{m}=
 \bbb{m}'\bbb{S}^{-1}\bbb{m}/(1+\bbb{m}'\bbb{S}^{-1}\bbb{m})$. \\
 {\rm (ii)} If either $n=2$ or the {\rm (free of parameters)} d.f.\
 $F(x)$ is such that
 \be \label{4.3}
 \Cov[Z_i,Z_j] \leq 0\ \ \mbox{
 for all $i\neq j$,\ \ $i,j=1,\ldots,n-1$,}
 \ee
 then the partial
 maxima {\rm BLUE (and BLIE)} of $\theta_2$ is non-negative.
 \end{theo}
 \vspace*{-0.6em}

 Note that, as in order statistics,
 the non-negativity of $L_2$
 is equivalent to the fact
 that the vector $\bbb{S}^{-1}\bbb{m}$
 (or, equivalently, the vector $\bbb{D}^{-1}\bbb{m}$) has
 non-negative entries; see Balakrishnan and Papadatos (2002) and Sarkadi (1985).

 Since it is important to know whether
 (\ref{4.3}) holds, in the sequel we shall make use of the following
 definition.

 \begin{DEFI}\label{def4.1}
 A d.f.\ $F(x)$ with finite second moment (or the
 corresponding density $f(x)$,
 if exists)
 \\
 (i) belongs to the class NCS (negatively correlated spacings)
 if its order statistics have negatively correlated
 spacings for all sample sizes $n\geq 2$.
 \\
 (ii)
 belongs
 to the class NCP if it has negatively correlated partial maxima
 spacings, i.e., if (\ref{4.3}) holds for all $n\geq 2$.
 \end{DEFI}

 An important result by Bai, Sarkar \& Wang (1997)
 states that a
 log-concave density $f(x)$ with finite variance belongs to NCS -- cf.\
 Sarkadi (1985). We
 call this sufficient condition as the S/BSW-condition
 (for ordinary spacings). Burkschat (2009, Theorem 3.5) showed an extended
 S/BSW-condition,
 under which the log-concavity of both $F$ and $1-F$ suffice
 for the NCS class.
 Due to the existence of simple
 formulae like (\ref{4.1}) and (\ref{4.2}), the NCS and NCP classes
 provide useful tools in verifying consistency for the scale
 estimator, as well as, non-negativity.
 Our purpose is to prove an S/BSW-type
 condition for partial maxima (see Theorem
 \ref{theo4.2}, Corollary \ref{cor4.1}, below). To this end,
 we first state
 Lemma \ref{lem4.1},
 that will be used in the sequel.
 %

 Only through the rest of the present section, we shall use
 the
 notation $Y_k=\max\{X_1,\ldots,X_k\}$,
 for any integer $k\geq 1$.

 \begin{lem}\label{lem4.1} Fix two integers $i$, $j$, with $1\leq i<j$, and suppose
 that the i.i.d.\ r.v.'s $X_1,X_2,\ldots$ have a common d.f.\
 $F(x)$.
 Let $I(\mbox{expression})$ denoting the indicator function taking
 the value $1$, if the expression holds true, and $0$ otherwise.
 \\ {\rm (i)} The conditional d.f.\ of $Y_{j+1}$ given
 $Y_j$ is
 \[
 \Pr[Y_{j+1} \leq y |\ Y_j]=\left\{
 \begin{array}{ccc}
 0, & \mbox{if} & y<Y_j\\
 F(y), & \mbox{if} & y\geq Y_j
 \end{array}
 \right.=F(y)I(y\geq Y_j),\ \ y\in\R.
 \]
 If, in addition, $i+1<j$, then the following property {\rm (which
 is an immediate consequence of the Markovian character of the
 extremal process)} holds:
 \[
 \Pr[Y_{j+1}\leq y |\ Y_{i+1}, Y_j]=\Pr[Y_{j+1}\leq y |\
 Y_j],\ \ y\in\R.
 \]
 {\rm (ii)} The conditional d.f.\ of $Y_i$ given
 $Y_{i+1}$ is
 \begin{eqnarray*}
 \Pr[Y_i \leq x |\ Y_{i+1}]
  & = & \left\{
 \begin{array}{ccc}
 \displaystyle \frac{F^i(x)}{\sum_{j=0}^i F^j(Y_{i+1})F^{i-j}(Y_{i+1}-)}, & \mbox{if} & x<Y_{i+1} \\
 1, & \mbox{if} & x\geq Y_{i+1}
 \end{array}
 \right. \\
 & = & I(x\geq Y_{i+1})
 +I(x<Y_{i+1})\frac{F^i(x)}{\sum_{j=0}^i F^j(Y_{i+1})F^{i-j}(Y_{i+1}-)},\ \ x\in\R.
 \end{eqnarray*}
 If, in addition, $i+1<j$, then the following
 property {\rm (which is again an immediate consequence of the
 Markovian character of the extremal process)} holds:
 \[
 \Pr[Y_i\leq x |\ Y_{i+1}, Y_j]=\Pr[Y_i \leq x |\
 Y_{i+1}],\ \ x\in\R.
 \]
 {\rm (iii)} Given $(Y_{i+1}, Y_{j})$, the random variables $Y_i$
 and $Y_{j+1}$ are independent.
 \end{lem}

 We omit the proof since the assertions are simple by-products of the
 Markovian character of the process $\{Y_k,\ k\geq 1\}$, which can
 be embedded in a continuous time extremal process $\{Y(t),\ t>0\}$;
 see Resnick (1987, Chapter 4).
 We merely note that a version of the Radon-Nikodym derivative
 of $F^{i+1}$ w.r.t.\ $F$ is given by
 \be
 \label{(4.4)}
 h_{i+1}(x)=\frac{dF^{i+1}(x)}{dF(x)}=\sum_{j=0}^{i}
 F^j(x)F^{i-j}(x-), \ \ \ x\in\R,
 \ee
 which is equal to $(i+1)F^{i}(x)$ only if  $x$ is a continuity point of
 $F$. To see this, it suffices to verify the identity
 \be
 \label{(4.5)}
 \int_{B}  \ dF^{i+1}(x)=\int_{B}  h_{i+1}(x)
 \ dF(x)  \ \
 \mbox{for all Borel sets} \ \ B\subseteq \R.
 \ee
 Now (4.5) is proved as follows: 
 \begin{eqnarray*}
 \int_{B}  \ dF^{i+1}  & = & \Pr(Y_{i+1}\in B) \\
  & = & \sum_{j=1}^{i+1} \Pr\left[ Y_{i+1}\in B, \ \sum_{k=1}^{i+1}
  I(X_k=Y_{i+1})=j\right] \\
  &=& \sum_{j=1}^{i+1} \sum_{1\leq k_1<\cdots<k_j\leq i+1}
  \Pr[ X_{k_1}=\cdots=X_{k_j}\in B, \\
  && \hspace{22ex}X_s<X_{k_1} \mbox{ for } s\notin
  \{k_1,\ldots,k_j\}] \\
  &=& \sum_{j=1}^{i+1} {i+1 \choose j}
  \Pr[ X_{1}=\cdots=X_j\in B, X_{j+1}<X_{1},\ldots,X_{i+1}<X_1] \\
  &=& \sum_{j=1}^{i+1} {i+1 \choose j}
  \int_B \E\left[\left(\prod_{k=2}^j I(X_k=x)\right)
  \left(\prod_{k=j+1}^{i+1} I(X_k<x)\right)\right] \ dF(x) \\
  &=&
  \int_B \left(\sum_{j=1}^{i+1} {i+1 \choose j} (F(x)-F(x-))^{j-1} F^{i+1-j}(x-) \right)
  \ dF(x) \\
  & = & \int_{B}  h_{i+1}(x) dF(x),
  \end{eqnarray*}
 where we used the identity $\sum_{j=1}^{i+1}{i+1 \choose j}(b-a)^{j-1} a^{i+1-j}=
 \sum_{j=0}^{i}b^{j} a^{i-j}$, $a\leq b$.

 We can now show the main result of this section, which presents an
 S/BSW-type condition for the partial maxima spacings.

 \begin{theo}\label{theo4.2}
 Assume that the
 d.f.\ $F(x)$, with
 finite second moment,
 is a log-concave distribution
 {\rm (in the sense that $\log F(x)$ is a concave
 function in $J$, where
 $J=\{x\in\R: 0<F(x)<1\}$)}, and has not an atom at its right
 end-point, $\omega(F)=\inf\{x\in\R :F(x)=1\}$. Then, $F(x)$ belongs to
 the class {\rm NCP}, i.e., {\rm (\ref{4.3})} holds for all $n\geq 2$.
 %
 \end{theo}
 \begin{pr}{Proof}
 For arbitrary r.v.'s $X\geq x_0>-\infty$ and $Y\leq y_0<+\infty$, with
 respective d.f.'s $F_X$, $F_Y$,
 we have
 \be\label{4.9}
 \E[X]=x_0+\int_{x_0}^{\infty} (1-F_X(t)) \ dt \ \
 \mbox{and}\ \ \ \E[Y]=y_0-\int_{-\infty}^{y_0} F_Y(t) \ dt
 \ee
 (cf.\ Papadatos (2001), Jones and Balakrishnan (2002)). Assume
 that $i<j$. By Lemma \ref{lem4.1}(i) and (\ref{4.9}) applied to
 $F_X=F_{Y_{j+1}| Y_j}$,
 it follows that
 \[
 \E [ Y_{j+1}|\ Y_{i+1},Y_j]=\E [ Y_{j+1}|\ Y_j]
 =Y_j+\int_{Y_j}^{\infty} (1-F(t)) \
 dt, \ \ \ \mbox{w.p.\ 1}.
 \]
 Similarly, by Lemma \ref{lem4.1}(ii) and (\ref{4.9}) applied to
 $F_Y=F_{Y_i| Y_{i+1}}$,
 we conclude that
 \[
 \E [ Y_i | \ Y_{i+1},Y_j]=\E [ Y_i | \ Y_{i+1}]
 =Y_{i+1}-\frac{1}{h_{i+1}(Y_{i+1})}
 \int_{-\infty}^{Y_{i+1}} F^i (t) \ dt, \ \ \ \mbox{w.p.\ 1},
 \]
 where $h_{i+1}$ is given by (\ref{(4.4)}). Note that $F$
 is continuous on $J$, since it is log-concave there, and thus,
 $h_{i+1}(x)=(i+1)F^i(x)$ for $x\in J$. If $\omega(F)$ is finite,
 $F(x)$ is also continuous at $x=\omega(F)$, by assumption.
 On the other hand, if $\alpha(F)=\inf\{x:F(x)>0\}$ is finite, $F$ can be
 discontinuous at
 $x=\alpha(F)$, but in this case, $h_{i+1}(\alpha(F))= F^i(\alpha(F))>0$; see
 (\ref{(4.4)}). Thus, in all cases, $h_{i+1}(Y_{i+1})>0$ w.p.\ 1.

 By conditional independence of $Y_i$ and $Y_{j+1}$ (Lemma
 \ref{lem4.1}(iii)), we have
 \begin{eqnarray*}
 \Cov(Z_i,Z_j|\ Y_{i+1},Y_j) & = &
 \Cov(Y_{i+1}-Y_i,Y_{j+1}-Y_j|\ Y_{i+1},Y_j) \\
 &=& -\Cov (Y_i,Y_{j+1}| Y_{i+1},Y_j)=0, \ \ \ \mbox{w.p.\ 1},
 \end{eqnarray*}
 so that $\E[\Cov(Z_i,Z_j|\ Y_{i+1},Y_j)]=0$, and thus,
 \begin{eqnarray}
 \Cov [Z_i,Z_j] & = &
 \Cov [\E(Z_i|\ Y_{i+1},Y_j), \E(Z_j|\ Y_{i+1},Y_j)]
 +\E [\Cov(Z_i,Z_j|\ Y_{i+1},Y_j)]
 \nonumber \\
 & = & \Cov [\E(Y_{i+1}-Y_i |\ Y_{i+1},Y_j), \E(Y_{j+1}-Y_j|\
    Y_{i+1},Y_j)]
 \nonumber \\
 & = & \Cov [Y_{i+1}-\E (Y_i |\ Y_{i+1},Y_j),
 \E(Y_{j+1}|\ Y_{i+1},Y_j)-Y_j]
 \nonumber \\
 &= & \Cov [ g(Y_{i+1}), h(Y_j) ],
 \label{(4.8)}
 \end{eqnarray}
 where
 \[
 g(x)=
 \left\{
 \begin{array}{ll}
 \displaystyle
 \frac{1}{(i+1)F^{i}(x)}\int_{-\infty}^x F^i(t) \ dt, & x>\alpha(F),
 \\
 \vspace{-1em}
 \\
 0, & \mbox{otherwise,}
 \end{array}\right.
 \ \ \
 h(x)=\int_{x}^{\infty} (1-F(t)) \ dt.
 \]
 Obviously, $h(x)$ is non-increasing. On the
 other hand, $g(x)$ is non-decreasing in $\R$.
 This can be
 shown as follows. First observe that $g(\alpha(F))=0$ if
 $\alpha(F)$ is finite, while $g(x)>0$ for $x>\alpha(F)$.
 Next observe that $g$ is finite and continuous at $x=\omega(F)$
 if $\omega(F)$ is finite, as follows by the assumed continuity of $F$ at
 $x=\omega(F)$ and the fact that $F$ has finite
 variance.
 Finally, observe that $F^i(x)$, a product of
 log-concave functions, is also log-concave in $J$. Therefore,
 for arbitrary $y\in J$, the function
 $d(x)=F^i(x)/\int_{-\infty}^y F^i(t)  dt$, $x\in(-\infty,y)\cap
 J$, is a probability density, and thus, it is a log-concave
 density with support $(-\infty,y)\cap J$.
 By Pr\`{e}kopa (1973) or Dasgupta and
 Sarkar (1982) it follows that the corresponding distribution
 function, $D(x)=\int_{-\infty}^x d(t) dt=\int_{-\infty}^x F^i(t)
 dt/ \int_{-\infty}^y F^i(t) dt$, $x\in (-\infty,y)\cap J$, is a
 log-concave distribution, and since $y$ is arbitrary,
 $H(x)=\int_{-\infty}^x F^i(t) dt$ is a log-concave function, for
 $x\in J$. Since $F$ is continuous in $J$, this is equivalent
 to the fact that the function
 \[
 \frac{H'(x)}{H(x)}=\frac{F^i(x)}{\int_{-\infty}^x F^i(t) \ dt}, \ \
 x\in J,
 \]
 is non-increasing, so that $g(x)=H(x)/((i+1)H'(x))$
 is non-decreasing in $J$.

 The desired result follows from (\ref{(4.8)}), because the r.v.'s
 $Y_{i+1}$ and $Y_j$ are positively quadrant dependent (PQD --
 Lehmann (1966)), since it is readily verified that $F_{Y_{i+1},Y_j}(x,y)\geq
 F_{Y_{i+1}}(x)F_{Y_j}(y)$ for all $x$ and $y$ (Lemma
 \ref{lem2.2}(i)). This completes the proof. $\Box$
 \end{pr}
 \bigskip

 The restriction $F(x)\to 1$ as $x\to \omega(F)$ cannot be removed
 from the theorem. Indeed, the function
 \[
 F(x)=\left\{
 \begin{array}{ll}
 0 & x\leq 0, \\
 x/4, & 0\leq x<1, \\
 1, & x\geq 1,
 \end{array}
 \right.
 \]
 is a log-concave distribution in $J=(\alpha(F),\omega(F))=(0,1)$,
 for which $\Cov[Z_1,Z_2]=\frac{59}{184320}>0$. The function $g$,
 used in the proof,
 is given by
 \[
 g(x)=
 \left\{
 \begin{array}{ll}
 \displaystyle
 \max\left\{0,\frac{x}{(i+1)^2}\right\}, &  x<1,
 \\
 \vspace{-1em}
 \\
 \displaystyle
 \frac{x-1}{i+1}+\frac{1}{(i+1)^2 4^i}, &  x\geq 1,
 \end{array}\right.
 \]
 and it is not monotonic.

 Since the family of densities with log-concave distributions
 contains both families of log-concave and non-increasing densities
 (see, e.g., Pr\`{e}kopa (1973), Dasgupta and Sarkar (1982), Sengupta and Nanda (1999),
 Bagnoli and Bergstrom (2005)), the following corollary is an immediate
 consequence of Theorems \ref{theo4.1} and \ref{theo4.2}.

 \begin{cor}\label{cor4.1}
 Assume that $F(x)$ has finite second moment.
 \\ {\rm (i)} If
 $F(x)$ is a log-concave
 d.f.\ {\rm (in
 particular, if $F(x)$ has either a log-concave or a non-increasing
 (in its interval support) density $f(x)$)},
 then the partial maxima {\rm BLUE} and the
 partial maxima {\rm BLIE} of $\theta_2$ are non-negative.
 \\ {\rm (ii)} If $F(x)$ has either a log-concave or a non-increasing
 {\rm (in its interval support)}
 density $f(x)$ then it belongs to the {\rm NCP} class.
 \end{cor}

 Sometimes it is asserted that ``the distribution of a log-convex density
 is log-concave'' (see, e.g., Sengupta and Nanda (1999), Proposition 1(e)),
 but this is not correct in its full generality, even if
 the corresponding r.v.\ $X$ is non-negative. For example, let
 $Y\sim$ Weibull with shape parameter $1/2$, and set $X\law Y|Y<1$.
 Then $X$ has density
 $f$ and d.f.\ $F$ given by
 \[
 f(x)=\frac{\exp(-\sqrt{1-x})}{2(1-e^{-1})\sqrt{1-x}}, \ \ \
 F(x)=\frac{\exp(-\sqrt{1-x})-e^{-1}}{1-e^{-1}}, \ \ \ 0<x<1,
 \]
 and it is easily checked that $\log f$ is convex in $J=(0,1)$, while
 $F$ is not log-concave in $J$. However, we point out that
 if $\sup J=+\infty$ then any log-convex density, supported on $J$, has to be
 non-increasing in $J$ and, therefore, its distribution is
 log-concave in $J$. Examples of log-convex distributions
 having a log-convex density are given by Bagnoli and Bergstrom (2005).



 \section{Consistent estimation of the scale parameter}\vspace*{-.5em}
 \setcounter{equation}{0}\label{sec5}

 Through this section we always assume that $F(x)$, the d.f.\ that
 generates the location-scale family, is non-degenerate and has
 finite second moment. The main purpose is to verify
 consistency for $L_2$, applying the results of section \ref{sec4}.
 To this end, we firstly state and prove a simple lemma that goes
 through the lines of Lemma \ref{lem2.1}. Due to the obvious fact
 that $\MSE[T_2]\leq \Var[L_2]$, all the results of the present
 section apply also to the BLIE of $\theta_2$.

 \begin{lem} \label{lem5.1}
 If $F(x)$ belongs to the {\rm NCP} class then \\
 {\rm (i)}
  \be\label{5.1}
 \Var [L_2]\leq \frac{\theta_2^2}{\sum_{k=1}^{n-1}
 m_k^2/s_k^2},
  \ee
 where $m_k=\E[Z_k]$ is the $k$-th component of the vector
 $\bbb{m}$ and $s_k^2=s_{kk}=\Var[Z_k]$ is the $k$-th diagonal
 entry of the matrix $\bbb{S}$. \\ {\rm (ii)} The partial maxima
 {\rm BLUE}, $L_2$, is consistent if the series
 \be \label{5.2}
 \sum_{k=1}^{\infty} \frac{m_k^2}{s_k^2}=+\infty.
 \ee
 \end{lem}
 \begin{pr}{Proof}
 Observe that part (ii) is an immediate consequence
 of part (i), due to the fact that, in contrast
 to the order statistics setup, $m_k$ and
 $s_k^2$ do not depend on the sample size $n$.
 Regarding (i), consider the linear
 unbiased estimator
 \[
 U_2=\frac{1}{c_n}\sum_{k=1}^{n-1} \frac{m_k}{s_k^2} Z_k^*\law
 \frac{\theta_2}{c_n} \sum_{k=1}^{n-1} \frac{m_k}{s_k^2} Z_k,
 \]
 where $c_n=\sum_{k=1}^{n-1} m_k^2/s_k^2$. Since $F(x)$ belongs to
 NCP and the weights of $U_2$ are positive, it follows that the
 variance of $U_2$, which is greater than or equal to the variance
 of $L_2$, is bounded by the RHS of (\ref{5.1}); this completes the
 proof. $\Box$
 \end{pr}
 \bigskip

 The proof of the following theorem is now immediate.
 \begin{theo}\label{theo5.1}
 If $F(x)$ belongs to the {\rm NCP}
 class and if there exists a finite constant $C$
 and a positive integer $k_0$
 such that
 \be\label{5.3}
 \frac{ \E[Z_k^2]}{k\E^2 [Z_k]}\leq C ,\  \mbox{ for
 all } k\geq k_0,
 \ee
 then
 \be\label{5.4}
 \Var [ L_2] \leq O\left( \frac{1}{\log n} \right),
 \ \mbox{ as }n\to\infty.
 \ee
 \end{theo}
 \begin{pr}{Proof}
 Since for $k\geq k_0$,
 \[
 \displaystyle
 \frac{m_k^2}{s_k^2}=\frac{m_k^2}{\E [Z_k^2]-m_k^2}=
 \frac{1}{\displaystyle
 \frac{\E[Z_k^2]}{\E^2[Z_k]}-1}\geq \frac{1}{Ck-1},
 \]
 the result follows by (\ref{5.1}). $\Box$
 \end{pr}
 \bigskip

 Thus, for proving consistency of order $1/\log n$ into NCP class
 it is sufficient to verify (\ref{5.3}) and, therefore, we shall
 investigate the quantities $m_k=\E[Z_k]$ and
 $\E[Z_k^2]=s_k^2+m_k^2$. A simple application of Lemma
 \ref{lem2.2}, observing that $m_k=\mu_{k+1}-\mu_k$ and
 $s_k^2=\sigma_{k+1}^2-2\sigma_{k,k+1}+\sigma_k^2$, shows that
 \begin{eqnarray}
 \label{5.5}\E[Z_k] & = & \int_{-\infty}^{\infty} F^{k}(x) (1-F(x)) \ dx, \\
 \label{5.6}\E^2[Z_k] & = & 2\int\int_{-\infty<x<y<\infty}
 F^{k}(x) (1-F(x)) F^{k}(y) (1-F(y)) \ dy \ dx,\\
 \label{5.7}\E[Z_k^2] & = & 2\int\int_{-\infty<x<y<\infty} F^k(x) (1-F(y))
 \ dy \
 dx.
 \end{eqnarray}
 Therefore, all the quantities of interest can be expressed as
 integrals in terms of the (completely arbitrary) d.f.\ $F(x)$
 (cf.\ Jones and Balakrishnan (2002)).

 For the proof of the main result we finally need the following
 lemma and its corollary.

 \begin{lem} \label{lem5.2}{\rm (i)} For any $t>-1$,
 \be\label{5.8}
 \lim_{k\to\infty} k^{1+t} \int_0^1 u^k (1-u)^t \
 du=\Gamma(1+t)>0.
 \ee
 {\rm (ii)}
 For any $t$ with $0\leq t<1$ and any $a>0$, there exist positive
 constants $C_1$, $C_2$, and a positive integer $k_0$ such that
 \be\label{5.9}
 0<C_1< k^{1+t} (\log k)^a \int_0^1 \frac{u^k (1-u)^t}{L^a(u)} \
 du<C_2<\infty, \ \mbox{ for all }k\geq k_0,
 \ee
 where $L(u)=-\log
 (1-u)$.
 \end{lem}
 \begin{pr}{Proof}
 Part (i) follows by Stirling's formula.
 For part (ii),
 with the substitution $u=1-e^{-x}$, we write the integral in
(\ref{5.9}) as
 \[
 \frac{1}{k+1}\int_0^{\infty} (k+1)(1-e^{-x})^k e^{-x}
 \frac{\exp(-tx)}{x^a} \ dx = \frac{1}{k+1}
 \E \left[ \frac{\exp(-t T)}{T^a}\right],
 \]
 where $T$ has the same distribution as the maximum of $k+1$ i.i.d.\
 standard exponential r.v.'s. It is well-known that
 $\E [T]=1^{-1}+2^{-1}+\cdots+(k+1)^{-1}$. Since the second derivative
 of the function $x\to e^{-tx}/x^a$ is $x^{-a-2}e^{-tx}
 (a+(a+tx)^2)$, which is positive for $x>0$, this function is
 convex, so by Jensen's inequality we conclude that
 \begin{eqnarray*}
 k^{1+t} (\log k)^a \int_0^1 \frac{u^k (1-u)^t}{(L(u))^a} \ du \geq
 \frac{k}{k+1} \left(
 \frac{\log k}{1+1/2+\cdots+1/(k+1)}\right)^a\hspace*{7ex}\mbox{~}
 \\
 \times \exp\left[-t\left(1+\frac12+\cdots+\frac{1}{k+1}-\log k \right)
 \right],
 \end{eqnarray*}
 and the RHS remains positive as $k\to\infty$, since it converges
 to $e^{-\gamma t}$, where $\gamma=.5772\ldots$ is Euler's
 constant. This proves the lower bound in (\ref{5.9}). Regarding
 the upper bound, observe that the function $g(u)=(1-u)^t/L^a(u)$,
 $u\in(0,1)$,
 has second derivative
 \[
 g''(u)=\frac{-(1-u)^{t-2}}{L^{a+2}(u)} [ t(1-t)
 L^2(u)+a(1-2t)L(u)-a(a+1)], \ \ 0<u<1,
 \]
 and since $0\leq t<1$, $a>0$ and $L(u)\to +\infty$ as $u\to 1-$,
 it follows that there exists a constant $b\in (0,1)$ such that
 $g(u)$ is concave in $(b,1)$. Split now the integral in
 (\ref{5.9}) in two parts,
 \[
 I_k=I_k^{(1)}(b)+I_k^{(2)}(b)=\int_0^b \frac{u^k (1-u)^t}{L^a(u)}
 \ du+\int_b^1 \frac{u^k (1-u)^t}{L^a(u)} \ du,
 \]
 and observe that for any fixed $s>0$ and any fixed
 $b\in(0,1)$,
 \[
 k^s I_k^{(1)}(b)\leq k^s b^{k-a} \int_0^b \frac{u^a
 (1-u)^t}{L^a(u)} \ du\leq k^s b^{k-a} \int_0^1 \frac{u^a
 (1-u)^t}{L^a(u)} \ du\to 0, \ \ \mbox{ as } k\to\infty,
 \]
 because the last integral is finite and independent of $k$.
 Therefore, $k^{1+t}(\log k)^a I_k$ is bounded above if
 $k^{1+t}(\log k)^a I_k^{(2)}(b)$ is bounded above for some $b<1$.
 Choose  $b$ close enough to $1$ so that $g(u)$ is concave in
 $(b,1)$.
 By Jensen's inequality and the fact that
 $1-b^{k+1}<1$ we conclude that
 \[
 I_k^{(2)}(b)=\frac{1-b^{k+1}}{k+1}\int_b^1 f_k(u) g(u) \ du
 =\frac{1-b^{k+1}}{k+1}\E[g(V)]\leq \frac{1}{k+1} g[\E(V)],
 \]
 where $V$ is an r.v.\ with density $f_k(u)=(k+1)u^k/(1-b^{k+1})$,
 for $u\in(b,1)$. Since
 $\E(V)=((k+1)/(k+2))(1-b^{k+2})/(1-b^{k+1})>(k+1)/(k+2)$, and $g$
 is positive and decreasing (its first derivative is
 $g'(u)=-(1-u)^{t-1}L^{-a-1}(u)(tL(u)+a)<0$, $0<u<1$), it follows
 from the above inequality that
 \begin{eqnarray*}
 k^{1+t}(\log k)^a I_k^{(2)}(b)\leq \frac{k^{1+t}(\log k)^a} {k+1}\
 g[\E(V)] \leq \frac{k^{1+t}(\log k)^a}{k+1}\
 g\left(\frac{k+1}{k+2}\right)
 \\
 =\frac{k}{k+1} \left( \frac{\log k }{\log (k+2)}\right)^a
 \left(\frac{k}{k+2}\right)^t \leq 1.
 \end{eqnarray*}
 This shows that $k^{1+t}(\log k)^a I_k^{(2)}(b)$ is bounded above,
 and thus, $k^{1+t}(\log k)^a I_k$ is bounded above,
 as it was to be shown. The proof
 is complete. $\Box$
 \end{pr}
 \begin{cor}\label{cor5.1} {\rm (i)} Under the assumptions of Lemma
 {\rm \ref{lem5.2}(i)}, for any $b\in [0,1)$ there exist positive
 constants $A_1$, $A_2$ and a positive integer $k_0$ such that
 \[
 0<A_1< k^{1+t} \int_b^1 u^k (1-u)^t \ du<A_2<+\infty,
 \ \mbox{ for all }k\geq k_0.
 \]
 {\rm (ii)} Under the assumptions of Lemma {\rm \ref{lem5.2}(ii)}, for
 any $b\in[0,1)$
 there exist positive
 constants $A_3$, $A_4$ and a positive integer $k_0$ such that
 \[
 0<A_3< k^{1+t} (\log k)^a \int_b^1 \frac{u^k (1-u)^t}{L^a(u)} \
 du<A_4<\infty, \ \mbox{ for all }k\geq k_0.
 \]
 \end{cor}
 \begin{pr}{Proof} The proof follows from Lemma \ref{lem5.2} in
 a trivial way, since the corresponding integrals over $[0,b]$ are
 bounded above by a multiple of $b^{k-a}$, of the form $A b^{k-a}$,
 with $A<+\infty$ being independent of $k$. $\Box$
 \end{pr}
 \bigskip

 We can now state and prove our main result:
 \begin{theo}\label{theo5.2}
 Assume that $F(x)$ lies in {\rm NCP}, and let $\omega=\omega(F)$
 be the upper end-point of the support of $F$, i.e., $\omega=
 \inf\{x\in \R: F(x)=1\}$, where $\omega=+\infty$ if $F(x)<1$ for
 all $x$. Suppose that $\lim_{x\to\omega-}F(x)=1$, and that $F(x)$
 is differentiable in a left neighborhood $(M,\omega)$ of $\omega$,
 with derivative $f(x)=F'(x)$ for $x\in (M,\omega)$. For
 $\delta\in\R$ and $\gamma\in\R$, define the {\rm (generalized hazard
 rate)} function
 \be\label{5.10}
 L(x)=L(x;\delta,\gamma;F)=\frac{f(x)}{(1-F(x))^{\gamma}
 (-\log (1-F(x)))^\delta},\ \ \
 x\in (M,\omega),
 \ee
 and set
 \[
 L_*=L_{*}(\delta,\gamma;F)=\liminf_{x\to\omega-}
 L(x;\delta,\gamma,F),\ \ L^*=L^*(\delta,
 \gamma;F)=\limsup_{x\to\omega-} L(x;\delta,\gamma,F).
 \]
 If either {\rm (i)} for some $\gamma<3/2$ and $\delta=0$,
 \vspace*{.5ex}  or {\rm (ii)} for some $\delta>0$ and some
 $\gamma$ with $1/2<\gamma\leq 1$,
 \be\label{5.11}
 0<L_*(\delta,\gamma;F)\leq L^*(\delta,\gamma;F)<+\infty,
 \ee
 then
 the partial maxima {\rm BLUE} $L_2$ {\rm (given by (\ref{2.2}) or (\ref{4.1}))}
 of the scale parameter $\theta_2$ is consistent
 and, moreover, $\Var [L_2]\leq O(1/\log n)$.
 \end{theo}
 \begin{pr}{Proof}
 First observe that for large enough $x<\omega$,  (\ref{5.11})
 implies that
 $f(x)>(L_*/2)(1-F(x))^{\gamma}(-\log(1-F(x)))^{\delta}>0$, so that
 $F(x)$ is eventually strictly increasing and continuous. Moreover,
 the derivative $f(x)$ is necessarily finite
 since
 $f(x)<2 L^* (1-F(x))^{\gamma}(-\log(1-F(x)))^{\delta}$.
 The assumption $\lim_{x\to\omega-}F(x)=1$ now shows that
 $F^{-1}(u)$ is uniquely defined in a left neighborhood of $1$,
 that $F(F^{-1}(u))=u$ for $u$ close to $1$, and that
 $\lim_{u\to 1-}F^{-1}(u)=\omega$. This,
 in turn, implies that $F^{-1}(u)$ is
 differentiable for $u$ close to $1$, with (finite) derivative
 $(F^{-1}(u))'=1/f(F^{-1}(u))>0$. In view of Theorem \ref{theo5.1},
 it suffices to verify (\ref{5.3}), and thus we seek for an upper
 bound on $\E[Z_k^2]$ and for a lower bound on $\E[Z_k]$. Clearly,
 (\ref{5.3}) will be deduced if we shall verify that, under (i),
 there exist finite constants $C_3>0$, $C_4>0$ such that
 \be\label{5.12}
 k^{3-2\gamma}\E[Z_k^2]\leq C_3 \ \ \  \mbox{ and } \ \ \
 k^{2-\gamma}\E[Z_k]\geq C_4,
 \ee
 for all large enough $k$. Similarly, (\ref{5.3}) will be verified
 if we show that, under (ii), there exist finite constants $C_5>0$
 and $C_6>0$ such that
 \be\label{5.13}
 k^{3-2\gamma}(\log k)^{2\delta}\E[Z_k^2]\leq C_5 \ \ \ \mbox{ and } \ \ \
 k^{2-\gamma}(\log k)^{\delta} \E[Z_k]\geq C_6,
 \ee
 for all large enough $k$. Since the integrands in the integral
 expressions (\ref{5.5})-(\ref{5.7}) vanish if $x$ or $y$ lies
 outside the set $\{ x\in \R: 0<F(x)<1\}$,  we have the equivalent
 expressions
 \begin{eqnarray}
 \label{5.14}\E[Z_k] & = & \int_{\alpha}^{\omega} F^{k}(x) (1-F(x)) \ dx, \\
 \label{5.15}\E[Z_k^2] & = & 2\int\int_{\alpha<x<y<\omega} F^k(x)
 (1-F(y))
 \ dy \ dx,
 \end{eqnarray}
 where $\alpha$ (resp., $\omega$) is the lower (resp., the upper)
 end-point of the support of $F$. Obviously,
 for any fixed $M$
 with $\alpha<M<\omega$
 and any
 fixed $s>0$, we have, as in  the proof of Lemma \ref{lem5.2}(ii),
 that
 \begin{eqnarray*}
 \lim_{k\to\infty} k^s \int_{\alpha}^{M} F^{k}(x) (1-F(x)) \ dx=0, \\
 \lim_{k\to\infty} k^s \int_{\alpha}^M\int_x^{\omega} F^k(x) (1-F(y))
 \ dy \ dx=0,
 \end{eqnarray*}
 because $F(M)<1$ and both integrals (\ref{5.14}), (\ref{5.15}) are
 finite for $k=1$, by the assumption that the variance is finite
 ((\ref{5.15}) with $k=1$ just equals to the variance of $F$; see
 also (\ref{2.9}) with $i=1$). Therefore, in order to verify
 (\ref{5.12}) and (\ref{5.13}) for large enough $k$,
 it is sufficient to replace $\E[Z_k]$ and $\E[Z_k^2]$, in both
 formulae (\ref{5.12}), (\ref{5.13}), by the integrals
 $\int^{\omega}_{M} F^{k}(x) (1-F(x))  dx$ and
 $\int^{\omega}_M\int_x^{\omega} F^k(x) (1-F(y)) dy dx $,
 respectively, for an arbitrary (fixed) $M\in(\alpha,\omega)$.
 Fix now $M\in(\alpha,\omega)$ so large
 that $f(x)=F'(x)$ exists and it is finite and strictly positive
 for all $x\in(M,\omega)$,
 and make the transformation $F(x)=u$ in the first integral, and
 the transformation $(F(x),F(y))=(u,v)$ in the second one. Both
 transformations are now one-to-one and continuous, because both
 $F$ and $F^{-1}$ are differentiable in their respective intervals
 $(M,\omega)$ and $(F(M),1)$, and their derivatives are finite and
 positive. Since $F^{-1}(u)\to \omega$ as $u\to 1-$, it is easily
 seen that (\ref{5.12}) will be concluded if it can be shown that
 for some fixed $b<1$ (which can be chosen arbitrarily close to
 $1$),
 \begin{eqnarray}
 \label{5.16}
 k^{3-2\gamma}\int^{1}_b \frac{u^k}{f(F^{-1}(u))}
 \left( \int_u^{1} \frac{1-v}{f(F^{-1}(v))} \ dv \right)  du & \leq
 &
 C_3 \ \ \  \mbox{ and } \\
 \label{5.17}
 k^{2-\gamma}\int^{1}_{b} \frac{u^{k}
 (1-u)}{f(F^{-1}(u))} \ du & \geq & C_4,
 \end{eqnarray}
 holds for all large enough $k$.
 Similarly,
 (\ref{5.13}) will be deduced if
 it will be proved that
 for some fixed $b<1$
 (which can be chosen arbitrarily close to $1$),
 \begin{eqnarray}
 \label{5.18}
 k^{3-2\gamma}(\log k)^{2\delta}\int^{1}_b
 \frac{u^k}{f(F^{-1}(u))} \left( \int_u^{1}
 \frac{1-v}{f(F^{-1}(v))} \ dv \right)  du & \leq & C_5 \ \ \
 \mbox{ and } \\
 \label{5.19}
 k^{2-\gamma}(\log k)^{\delta} \int^{1}_{b} \frac{u^{k}
 (1-u)}{f(F^{-1}(u))} \ du & \geq & C_6,
 \end{eqnarray}
 holds for all large enough $k$. The rest of the proof is thus
 concentrated on showing (\ref{5.16}) and (\ref{5.17}) (resp.,
 ((\ref{5.18}) and (\ref{5.19})), under the assumption (i) (resp.,
 under the assumption (ii)).

 Assume first that (\ref{5.11}) holds under (i).
 Fix now $b<1$ so large that
 \[
 \frac{L_*}{2}(1-F(x))^{\gamma}<f(x) <2L^* (1-F(x))^{\gamma},\
 \mbox{ for all }x\in (F^{-1}(b),\omega);
 \]
 equivalently,
 \be\label{5.20} \frac{1}{2L^*}<
 \frac{(1-u)^{\gamma}}{f(F^{-1}(u))} <\frac{2}{L_*}, \ \ \mbox{ for
 all }u\in(b,1).
 \ee
 Due to (\ref{5.20}), the inner integral in (\ref{5.16}) is
 \[
  \int_u^{1}
 \frac{1-v}{f(F^{-1}(v))}
 \ dv  =
 \int_u^{1} (1-v)^{1-\gamma}
 \frac{(1-v)^{\gamma}}{f(F^{-1}(v))} \ dv \leq
 \frac{2(1-u)^{2-\gamma}}{(2-\gamma)L_*}.
 \]
 By Corollary \ref{cor5.1}(i) applied for $t=2-2\gamma>-1$, the LHS
 of (\ref{5.16}) is less than or equal to
 \[
 \frac{2k^{3-2\gamma}}{(2-\gamma)L_*} \int^{1}_b u^k
 (1-u)^{2-2\gamma}\frac{(1-u)^{\gamma}}{f(F^{-1}(u))} \ du \leq
 \frac{4k^{3-2\gamma}}{(2-\gamma)L_*^2} \int^{1}_b u^k
 (1-u)^{2-2\gamma} \ du \leq C_3,
 \]
 for all $k\geq k_0$, with $C_3=4 A_2
 L_*^{-2}(2-\gamma)^{-1}<\infty$, showing (\ref{5.16}). Similarly,
 using the lower bound in (\ref{5.20}), the integral in
 (\ref{5.17}) is
 \[
 \int^{1}_{b} \frac{u^{k} (1-u)}{f(F^{-1}(u))} \ du = \int^{1}_{b}
 u^k (1-u)^{1-\gamma}\frac{(1-u)^{\gamma}}{f(F^{-1}(u))}\ du \geq
 \frac{1}{2L^*} \int^{1}_{b} u^k (1-u)^{1-\gamma}\ du,
 \] so that, by Corollary \ref{cor5.1}(i) applied for
 $t=1-\gamma>-1$, the LHS of (\ref{5.17}) is greater than  or equal
 to
 \[
 \frac{k^{2-\gamma}}{2L^*} \int_b^1 u^k (1-u)^{1-\gamma} \ du \geq
 \frac{A_1}{2L^*}>0,\ \ \mbox{ for all $k\geq k_0$,}
 \] showing (\ref{5.17}).

  Assume now that (\ref{5.11}) is satisfied under (ii). As in part (i),
 choose a large enough $b<1$ so that
 \be\label{5.21}
 \frac{1}{2L^*}< \frac{(1-u)^{\gamma}L^{\delta} (u)}{f(F^{-1}(u))}
 <\frac{2}{L_*}, \ \ \mbox{ for all }u\in(b,1),
 \ee
 where $L(u)=-\log(1-u)$.
 Due to (\ref{5.21}), the inner integral in
 (\ref{5.18}) is
 \[
 \int_u^{1} \frac{(1-v)^{1-\gamma}}{L^{\delta}(v)}
 \frac{(1-v)^{\gamma}L^{\delta}(v)}{f(F^{-1}(v))} \ dv \leq
 \frac{2}{L_*}
 \int_u^{1} \frac{(1-v)^{1-\gamma}}{L^{\delta}(v)}
 \ dv
 \leq \frac{2(1-u)^{2-\gamma}}{L_* L^{\delta}(u)},
 \]
 because $(1-u)^{1-\gamma}/L^{\delta}(u)$ is decreasing (see the
 proof of Lemma \ref{lem5.2}(ii)). By Corollary \ref{cor5.1}(ii)
 applied for $t=2-2\gamma\in[0,1)$ and $a=2\delta>0$, the double
 integral in (\ref{5.18}) is less than or equal to
 \[
 \frac{2}{L_*} \int^{1}_b \frac{u^k
 (1-u)^{2-2\gamma}}{L^{2\delta}(u)} \frac{(1-u)^{\gamma}L^{\delta}
 (u)}{f(F^{-1}(u))}  du
 \leq \frac{4}{L_*^2} \int^{1}_b \frac{u^k
 (1-u)^{2-2\gamma}}{L^{2\delta}(u)} du \leq
 \frac{C_5}{k^{3-2\gamma}(\log k)^{2\delta}},
 \]
 for all $k\geq k_0$, with $C_5=4 A_4 L_*^{-2}<\infty$, showing
 (\ref{5.18}). Similarly, using the lower bound in (\ref{5.21}),
 the integral in (\ref{5.19}) is
 \[
 \int^{1}_{b} \frac{u^k (1-u)^{1-\gamma}}{L^{\delta}(u)}
 \frac{(1-u)^{\gamma}L^{\delta}(u)}{f(F^{-1}(u))}\ du \geq
  \frac{1}{2L^*} \int^{1}_{b} \frac{u^k
 (1-u)^{1-\gamma}}{L^{\delta}(u)}\ du,
 \]
 and thus, by Corollary \ref{cor5.1}(ii) applied for
 $t=1-\gamma\in[0,1)$ and $a=\delta>0$, the LHS of (\ref{5.19}) is
 greater than  or equal to
 \[
 \frac{k^{2-\gamma}(\log k)^{\delta}}{2L^*} \int_b^1 \frac{u^k
 (1-u)^{1-\gamma}}{L^{\delta}(u)} \ du \geq \frac{A_3}{2L^*}>0,\ \
 \mbox{ for all $k\geq k_0$,}
 \]
 showing (\ref{5.19}). This completes the proof. $\Box$
\end{pr}

 \begin{REM}\label{rem5.1}
 Taking $L(u)=-\log (1-u)$, the limits $L_*$ and $L^*$ in
 (\ref{5.11}) can be rewritten as
 \begin{eqnarray*}
 L_*(\delta,\gamma;F) & = & \liminf_{u\to
 1-}\frac{f(F^{-1}(u))}{(1-u)^{\gamma}L^{\delta}(u)}
 =\left(\limsup_{u\to
 1-}(F^{-1}(u))'(1-u)^{\gamma}L^{\delta}(u)\right)^{-1},
 \\
 L^*(\delta,\gamma;F)&=&\limsup_{u\to
 1-}\frac{f(F^{-1}(u))}{(1-u)^{\gamma}L^{\delta}(u)}
 =\left(\liminf_{u\to
 1-}(F^{-1}(u))'(1-u)^{\gamma}L^{\delta}(u)\right)^{-1}.
 \end{eqnarray*}
 In the particular case where $F$ is absolutely continuous with
 a continuous density $f$ and interval support, the function
 $f(F^{-1}(u))=1/(F^{-1}(u))'$ is
 known as the density-quantile function (Parzen, (1979)), and
 plays a fundamental role in the theory of order statistics.
 Theorem \ref{theo5.2} shows, in some sense, that the behavior of
 the density-quantile function at the upper end-point, $u=1$,
 specifies the variance behavior of the partial
 maxima BLUE for the scale parameter $\theta_2$. In fact,
 (\ref{5.11}) (and (\ref{6.1}), below) is a Von Mises-type
 condition (cf.\ Galambos (1978), \S\S 2.7, 2.11).
 \end{REM}
 \begin{REM}\label{rem5.2}
 It is obvious that condition $\lim_{x\to\omega-}F(x)=1$ is necessary for
 the consistency of BLUE (and BLIE). Indeed,
 the event that all partial maxima are equal to $\omega(F)$ has probability
 $p_0=F(\omega)-F(\omega-)$ (which is independent of $n$). Thus, a point mass at
 $x=\omega(F)$ implies that for all $n$, $\Pr(L_2=0)\geq p_0>0$. This situation is trivial.
 Non-trivial cases also exist, and we provide one at the end of next
 section.
 \end{REM}

 \section{Examples and conclusions}\vspace*{-.5em}
 \setcounter{equation}{0}\label{sec6}

 In most commonly used location-scale families, the following
 corollary suffices for concluding consistency of the BLUE (and the
 BLIE) of the scale parameter.
 Its proof follows by a
 straightforward combination of Corollary \ref{cor4.1}(ii) and
 Theorem \ref{theo5.2}.\vspace*{-0.6em}
\begin{cor}\label{cor6.1}
 Suppose that $F$ is absolutely continuous with finite variance,
 and that its density $f$ is
 either log-concave
 or non-increasing in its interval support
 $J=(\alpha(F),\omega(F))=\{ x\in\R: 0<F(x)<1\}$.
 If, either  for some $\gamma<3/2$ and $\delta=0$,
 or for some $\delta>0$ and some $\gamma$ with
 $1/2<\gamma\leq 1$,
 %
 \be\label{6.1}
 \lim_{x\to\omega(F)-} \frac{f(x)}{(1-F(x))^{\gamma} (-\log
 (1-F(x)))^\delta}=L\in(0,+\infty),
 \ee
 then the partial maxima
 {\rm BLUE} of the scale parameter is consistent and, moreover, its
 variance is at most of order $O(1/\log n)$.
 \end{cor}
 Corollary \ref{cor6.1} has immediate applications to several
 location-scale families. The following are some of them, where
 (\ref{6.1}) can be verified easily. In all these families
 generated by the distributions mentioned below, the variance of
 the partial maxima BLUE $L_2$ (see (\ref{2.3}) or (\ref{4.1})),
 and the mean squared error of the partial maxima BLIE $T_2$ (see
 (\ref{2.5}) or (\ref{4.2})) of the scale parameter is at most of
 order $O(1/\log n)$, as the sample size
 $n\to\infty$.\vspace*{-0.6em}
 \bigskip

 \noindent {\small\bf 1. Power distribution  (Uniform).}
 $F(x)=x^\lambda$, $f(x)=\lambda x^{\lambda-1}$, $0<x<1$
 ($\lambda>0$), and $\omega(F)=1$. The density is
 non-increasing for
 $\lambda\leq 1$ and log-concave for $\lambda\geq 1$. It is easily
 seen that (\ref{6.1}) is satisfied for $\delta=\gamma=0$ (for
 $\lambda=1$ (Uniform) see section \ref{sec3}).
 \bigskip

 \noindent {\small\bf 2. Logistic distribution.}
 $F(x)=(1+e^{-x})^{-1}$, $f(x)=e^{-x}(1+e^{-x})^{-2}$, $x\in\R$,
 and $\omega(F)=+\infty$. The density is  log-concave, and it is
 easily seen that (\ref{6.1}) is satisfied for $\delta=0$,
 $\gamma=1$.\vspace*{-0.6em}
 \bigskip

 \noindent {\small\bf 3. Pareto  distribution.} $F(x)=1-x^{-a}$,
 $f(x)=a x^{-a-1}$, $x>1$ ($a>2$, so that the second moment is
 finite), and $\omega(F)=+\infty$. The density is  decreasing, and
 it is easily seen that (\ref{6.1}) is satisfied for $\delta=0$,
 $\gamma=1+1/a$.
 Pareto case provides an example
 which lies in NCP and not in NCS class -- see
 Bai, Sarkar \& Wang (1997).
 \vspace*{-0.6em}
 \bigskip

 \noindent
 {\small\bf 4. Negative Exponential distribution.} $F(x)=f(x)=e^x$,
 $x<0$, and $\omega(F)=0$.
 The density is log-concave and it is
 easily seen that (\ref{6.1}) is satisfied for $\delta=\gamma=0$.
 This model is particularly important, because it corresponds to
 the partial minima model from the standard exponential
 distribution -- see Samaniego and Whitaker (1986).
 \bigskip

 \noindent
 {\small\bf 5. Weibull  distribution (Exponential).} $F(x)=1-e^{-x^c}$,
 $f(x)=cx^{c-1}\exp(-x^c)$,
 $x>0$ ($c>0$), and $\omega(F)=+\infty$.
 The density is  non-increasing for $c\leq 1$ and log-concave for
 $c\geq 1$, and it is easily seen that (\ref{6.1}) is satisfied for
 $\delta=1-1/c$,  $\gamma=1$. It should be noted that Theorem
 \ref{theo5.2} does not apply for $c<1$, since
 $\delta<0$.\vspace*{-0.6em}
 \bigskip

 \noindent
 {\small\bf 6. Gumbel (Extreme Value) distribution.}
 $F(x)=\exp(-e^{-x})=e^{x}f(x)$,
 $x\in\R$, and $\omega(F)=+\infty$.
 The distribution is log-concave and  (\ref{6.1})
 holds with $\gamma=1$, $\delta=0$ ($L=1$).
 This model is particularly important for its applications
 in forecasting records, especially in athletic events --
 see Tryfos and
 Blackmore
 \vspace{0.6em}
 (1985).

 \noindent {\small\bf 7. Normal Distribution.} $f(x)=\varphi(x)=(2\pi
 e^{x^2})^{-1/2}$, $F=\Phi$,
 $x\in \R$, and $\omega(F)=+\infty$.
 The density is log-concave and Corollary \ref{cor6.1} applies with
 $\delta=1/2$ and $\gamma=1$. Indeed,
 \[
 \lim_{+\infty} \frac{\varphi(x)}{(1-\Phi(x))(-\log(1-\Phi(x)))^{1/2}}=
 \lim_{+\infty} \frac{\varphi(x)}{x(1-\Phi(x))}\
 \frac{x}{(-\log(1-\Phi(x)))^{1/2}},
 \]
 and it is easily seen that
 \[
 \lim_{+\infty} \frac{\varphi(x)}{x(1-\Phi(x))}=1,\ \ \ \ \lim_{+\infty}
 \frac{x^2}{-\log(1-\Phi(x))}=2,
 \]
 so that $L=\sqrt{2}$.\vspace*{-0.6em}
 \bigskip

 In many cases of interest
 (especially in athletic events),
 best performances
 are presented as partial minima
 rather than maxima;
 see, e.g., Tryfos and Blackmore (1985).
 Obviously, the present theory applies also to
 the partial minima setup. The easiest way to convert
 the present results for the partial minima case is to consider
 the i.i.d.\ sample $-X_1,\ldots,-X_n$, arising from
 $F_{-X}(x)=1-F_X(-x-)$, and to observe that
 $ \min\{X_1,\ldots,X_i\}=-\max\{-X_1,\ldots,-X_i\}$, $i=1,\ldots,n$. Thus,
 we just have to replace $F(x)$ by $1-F(-x-)$
 in the corresponding formulae.

 There are some related problems and questions
 that, at least to our point of view, seem to be quite
 interesting.
 One problem is to verify consistency for the partial maxima BLUE
 of the location parameter. Another problem concerns the complete
 characterizations of the NCP and NCS classes (see Definition
 \ref{def4.1}), since we only know S/BSW-type sufficient conditions.
 Also, to prove or disprove the non-negativity of the partial maxima
 BLUE for the scale parameter, outside the NCP class (as well as
 for the order statistics BLUE of the scale parameter outside the
 NCS class).

 Some questions concern lower variance bounds for the
 partial maxima BLUEs. For example we showed in section \ref{sec3}
 that the rate $O(1/\log n)$ (which, by Theorem \ref{theo5.2},
 is just an upper bound for the variance of $L_2$)
 is the correct order for the
 variance of both estimators in the Uniform location-scale family.
 Is this the usual case?  If it is so, then we could properly
 standardize the estimators, centering and
 multiplying them by $(\log n)^{1/2}$.
 This would result to limit theorems analogues
 to the corresponding ones for order  statistics --
 e.g., Chernoff, Gastwirth \& Johns (1967); Stigler (1974);  --
 or analogues to the corresponding ones of Pyke (1965), (1980),
 for partial maxima spacings instead of ordinary spacings.
 However, note that the
 Fisher-Information approach, in the particular case of the one-parameter
 (scale) family generated by the standard Exponential distribution,
 suggests a variance of about $3\theta_2^2/(\log n)^3$ for the
 minimum variance unbiased estimator
(based on partial maxima) --
 see Hoffman and Nagaraja (2003, eq.\ (15) on p.\ 186).

 A final question concerns the
 construction of approximate BLUEs (for both location and scale)
 based on partial maxima,
 analogues to Gupta's (1952) simple linear estimators based on order
 statistics. Such a kind of approximations and/or limit theorems
 would be especially useful for practical purposes, since the
 computation of BLUE via its closed formula requires inverting
 an $n\times n$ matrix. This problem has been partially
 solved here: For the NCP class, the
 estimator $U_2$, given in the proof of Lemma \ref{lem5.1}, is
 consistent for $\theta_2$ (under the assumptions of Theorem
 \ref{theo5.2}) and can be computed by a simple formula if we
 merely know the means and variances of the partial maxima
 spacings.

 Except of the trivial case given in Remark \ref{rem5.2}, above, there exist non-trivial
 examples where no consistent sequence of unbiased estimators exist for the scale parameter.
 To see this, we make use of the following result.

 \begin{theo}\label{theo6.1} {\rm (Hofmann and Nagaraja 2003, p.\
 183)} Let $\Xs$ be an i.i.d.\ sample
 from the scale family with distribution function
 $F(x;\theta_2)=F(x/\theta_2)$ {\rm ($\theta_2>0$ is the scale parameter)} and density
 $f(x;\theta_2)=f(x/\theta_2)/\theta_2$, where $f(x)$ is known,
 it has a continuous derivative $f'$, and its support, $J(F)=\{x:f(x)>0\}$,
 is one of the intervals $(-\infty,\infty)$, $(-\infty,0)$ or $(0,+\infty)$.
 \\
 {\rm (i)} The Fisher Information contained in the partial maxima
 data $\Xspo$ is given by
 \[
 I^{\max}_n=\frac{1}{\theta_2^2}\sum_{k=1}^n
 \int_{J(F)} f(x) F^{k-1}(x)\left(1+\frac{xf'(x)}{f(x)}+\frac{(k-1)xf(x)}{F(x)}\right)^2
 dx.
 \]
 {\rm (ii)} The Fisher Information contained in the partial minima
 data $\Xspmo$ is given by
 \[
 I^{\min}_n=\frac{1}{\theta_2^2}\sum_{k=1}^n
 \int_{J(F)} f(x) (1-F(x))^{k-1}\left(1+\frac{xf'(x)}{f(x)}-\frac{(k-1)xf(x)}{1-F(x)}\right)^2
 dx.
 \]
 \end{theo}

 It is clear that for fixed $\theta_2>0$, $I^{\max}_n$ and $I^{\min}_n$ both increase
 with the sample size $n$. In particular, if $J(F)=(0,\infty)$ then,
 by Beppo-Levi's Theorem,
 $I_n^{\min}$ converges (as $n\to\infty$) to
 its limit
 \be
 \label{eq6.2}
 I^{\min}=\frac{1}{\theta_2^2} \int_{0}^{\infty} \left\{\mu(x)
 \left(1+\frac{xf'(x)}{f(x)}-x\mu(x)\right)^2 +x^2
 \mu^2(x)\left(\lambda(x)+\mu(x)\right)\right\}
 dx,
 \ee
 where $\lambda(x)=f(x)/(1-F(x))$ and $\mu(x)=f(x)/F(x)$ is the
 failure rate and reverse failure rate of $f$, respectively.
 Obviously, if $I^{\min}<+\infty$, then the Cram\'{e}r-Rao
 inequality shows that
 no consistent sequence of unbiased estimators exists. This, of
 course, implies that in the corresponding scale family,
 any sequence of linear (in partial minima) unbiased estimators is inconsistent. The same is clearly true for the location-scale family, because any linear unbiased estimator for $\theta_2$ in the location-scale family is also a linear unbiased estimator for $\theta_2$ in the corresponding scale family.

 In the following we show that there exist distributions with
 finite variance such that $I^{\min}$ in (\ref{eq6.2}) is finite:
 Define $s=e^{-2}$ and
 \[
 F(x)=\left\{
 \begin{array}{cc}
 0, & x\leq 0, \vspace{.2em}\\
 \displaystyle
 \frac{1}{1-\log(x)}, & 0<x \leq s, \\
 & \vspace*{-.7em}\\
 1-(ax^2+bx+c) e^{-x}, & x\geq s,
 \end{array}
 \right.
 \]
 where
 \begin{eqnarray*}
 a&=&\frac{1}{54} \exp(e^{-2})(18-6e^2+e^4)\simeq 0.599, \\
 b&=&-\frac{2}{27} \exp(-2+e^{-2})(9-12e^2+2e^4)\simeq -0.339, \\
 c&=&\frac{1}{54} \exp(-4+e^{-2})(18-42e^2+43e^4)\simeq 0.798. \\
 \end{eqnarray*}
 Noting that $F(s)=1/3$, $F'(s)=e^2/9$ and $F''(s)=-e^4/27$,
 it can be easily verified that the corresponding density
 \[
 f(x)=\left\{
 \begin{array}{cc}
 \displaystyle
 \frac{1}{x(1-\log(x))^2}, & 0<x \leq s, \\
 & \vspace*{-.7em}\\
 (ax^2+(b-2a)x+c-b) e^{-x}, & x\geq s,
 \end{array}
 \right.
 \]
 is strictly positive for $x\in(0,\infty)$, processes
 finite moments of any order, and has continuous derivative
 \[
 f'(x)=\left\{
 \begin{array}{cc}
 \displaystyle
 \frac{1+\log(x)}{x^2(1-\log(x))^3}, & 0<x \leq s, \\
 & \vspace*{-.7em}\\
 -(ax^2+(b-4a)x+2a-2b+c) e^{-x}, & x\geq s.
 \end{array}
 \right.
 \]
 Now the integrand in (\ref{eq6.2}), say $S(x)$, can be written as
 \[
 S(x)=\left\{
 \begin{array}{cc}
 \displaystyle
 \frac{1-2\log(x)}{x(-\log(x))(1-\log(x))^3}, & 0<x \leq s, \\
 & \vspace*{-.7em}\\
 A(x)+B(x), & x\geq s,
  \end{array}
  \right.
 \]
 where, as $x\to+\infty$,
 $A(x)\sim A x^4 e^{-x}$ and $B(x)\sim B x^6 e^{-2x}$,  with $A$, $B$
 being positive constants. Therefore,
 $\int_{0}^s S(x)dx=\int_{2}^{\infty} \frac{1+2y}{y(1+y)^3}dy=\log(3/2)-5/18\simeq
 0.128$. Also, since $S(x)$ is continuous in $[s,+\infty)$ and
 $S(x)\sim A x^4e^{-x}$ as $x\to+\infty$, it follows that $\int_{s}^{\infty}
 S(x)dx<+\infty$ and $I^{\min}$ is finite.

 Numerical integration shows that $\int_{s}^{\infty} S(x)dx\simeq 2.77$
 and thus, $I^{\min}\simeq 2.9/\theta_2^2<3/\theta_2^2$.
 In view of the Cram\'{e}r-Rao bound
 this means that, even if a huge sample of partial minima has been recorded,
 it is impossible to construct an unbiased scale estimator
 with variance less than $\theta_2^2/3$. Also, it should be noted that a similar example
 can be constructed such that $f''(x)$ exists (and is continuous) for all $x>0$.
 
 Of course the above example can be adapted to the partial maxima case by considering the location-scale family generated by
 the distribution function
 \[
 F(x)=\left\{
 \begin{array}{cc}
 (ax^2-bx+c) e^{x}, & x\leq -s,
  \vspace*{.5em}
 \\
 \displaystyle
 \frac{-\log(-x)}{1-\log(-x)}, & -s\leq x <0, 
  \vspace*{.5em}
 \\
 1, & x\geq 0,
 \end{array}
 \right.
 \]
 with $s$, $a$, $b$ and $c$ as before.

 \vspace{-0.6em}
 \bigskip

 \noindent {\small\bf Acknowledgements.}
 Research partially
 supported by the University of Athens' Research Fund under Grant
 70/4/5637. Thanks are due to an anonymous referee for the very careful reading of the manuscript and also for correcting some mistakes.
 I would also like to express my thanks to Professors Fred
 Andrews, Mohammad Raqab, Narayanaswamy Balakrishnan, Barry Arnold
 and Michael Akritas, for their helpful discussions and comments;
 also to Professors Erhard Cramer and Claude Lef\`{e}vre, for
 bringing to my attention many interesting references related to
 log-concave densities. Special thanks are due to Professor Fredos
 Papangelou for his helpful remarks that led to a more general
 version of Theorem \ref{theo4.2}, and also to Dr.\ Antonis Economou,
 who's question for the Normal Distribution was found to be crucial
 for the final form of the main result.
 
 \end{document}